\documentclass[a4paper,reqno]{amsart}
\usepackage{amssymb}
\usepackage{amsmath}
\usepackage{amsfonts}
\usepackage{latexsym}

\theoremstyle{plain}
\newtheorem{theorem}{Theorem}
\newtheorem{definition}{Definition}
\numberwithin{equation}{section}
\input{tcilatex}

\begin{document}
\title[Quantum stochastic processes]{Reconstruction theorem for quantum
stochastic processes}
\author{V P Belavkin}
\address{Department of Applied Mathematics\\
MIEM, Moscow 109028, USSR}
\email{vpb@maths.nott.ac.uk}
\date{Submitted June 15, 1981, revision submitted March 13, 1984}
\thanks{This paper is translated from: \textit{Teoreticheskaya i
Matematicheskaya Fizika}, \textbf{62} (3) 409 -- 431 (1985)}
\keywords{Quantum stochastic process, Strong and weak equivalence,
Relativistic causal sets, Quantum covariant measurements, Quantum strong and
weak Markovianity.}

\begin{abstract}
Statistically interpretable axioms are formulated that define a quantum
stochastic process (QSP) as a causally ordered field in an arbitrary
space--time region $T$ of an open quantum system under a sequential
observation at a discrete space-time localization. It is shown that to every
QSP described in the weak sense by a self-consistent system of causally
ordered correlation kernels there corresponds a unique, up to unitary
equivalence, minimal QSP in the strong sense. It is shown that the proposed
QSP construction, which reduces in the case of the linearly ordered discrete 
$T=\mathbb{Z}$ to the construction of the inductive limit of Lindblad's
canonical representations \cite{8}, corresponds to Kolmogorov's classical
reconstruction \cite{12} if the order on $T$ is ignored and leads to Lewis
construction \cite{14} if one uses the system of all (not only causal)
correlation kernels, regarding this system as lexicographically preordered
on $\mathbb{Z}\times T$. The approach presented encompasses both
nonrelativistic and relativistic irreversible dynamics of open quantum
systems and fields satisfying the conditions of local commutativity
semigroup covariance. Also given are necessary and sufficient conditions of
dynamicity (or conditional Markovianity) and regularity, these leading to
the properties of complete mixing (relaxation) and ergodicity of the QSP.
\end{abstract}

\maketitle

\section{Introduction}

The problem of statistical foundation for the irreversible processes in open
systems of quantum thermodynamics and measurement theory encountered in
coherent optics, quantum communications and microelectronics \cite{1}--\cite%
{3} requires the development of a general theory of quantum stochastic
processes (QSP) that contains the classical theory and the known QSP models
as special cases. This theory must be operational and should admit a
microscopically consistent statistical interpretation for the irreversible
successive physical maps like quantum dynamical transformations and quantum
measurement operations. On the phenomenological level such quantum
operations were introduced already by von Neumann \cite{4} and considered in
a more general framework by Haag and Kastler \cite{5} and Davies and Lewis 
\cite{6}. A microscopically consistent operational approach to the general
QSP theory which will be followed here, containing the classical and quantum
statistical theories as special cases, was outlined in \cite{7}.

In physical applications, QSPs are usually described by special
chronologically ordered correlation functions, the only functions which can
be dynamically defined and tested on the basis of the statistics of
successive measurements in real time, called causal. An axiomatic definition
of QSPs based on causal correlation operators corresponding to a discrete
time $T=\mathbb{Z}$ was given by Lindblad \cite{8} for the case of a simple
observable algebra $\mathcal{B}=\mathcal{B}(\mathcal{H})$, and by author 
\cite{9} for a Boolean algebra $\mathcal{B}$ of measurable events $%
B\subseteq \mathrm{E}$. Such an approach leads to the description of
representations $\pi _{t},t\in T$ of the algebra $\mathcal{B}$ in an
expanding physical system, the mathematical definition of which lies in the
basis of the operational approach \cite{7} and will be developed here.

The aim of the present paper is threefold. First, we give a physically
interpretable --- in real time --- definition of a QSP as a family of
representations of the observable algebra $\mathcal{B}$ in a single (large)
quantum system by indicating a universal method of constructing
(reconstructing) such a system from \emph{causal correlation kernels}
described by the consistent axioms formulated in the paper. Second, we
reconstruct the unitary representation of the irreversible \emph{endomorphic
dynamics} corresponding to the stationary QSP and encompass in a unified
manner both nonrelativistic and relativistic covariant QSPs describing open
quantum systems and fields in a causally ordered space-time region $%
T\subseteq \mathbb{R}^{d+1}$ with respect to a given semigroup of symmetries
on $T$. And third, we derive the \emph{principle of nondisturbance} (or the
nondemolition principle) as microcausality of a given quantum subsystem with
respect to the successive measurements of a QSP, assumed to be accessible
for observations in the given space-time $T$.

In order to give an explicit statistical interpretation of a QSP as a
process of chronologically ordered measurements at arbitrary times $%
\{t_{1},\dots ,t_{n}\}\subset T$, we shall consider here only the QSP over
an \emph{event algebra}, that is a Boolean algebra $\mathcal{B}$, for which
the probabilities of finite sequences of events $B_{t_{1}},\dots
,B_{t_{n}}\in \mathcal{B}$, observed in `space--time" $t\in T$ are
determined directly by the diagonal values $\mu (B_{t_{1}},\dots ,B_{t_{n}})$
of the corresponding correlation kernels. It will be shown that the
self-consistent family of all causal correlation kernels that satisfy
covariance conditions that generalize Lorentz covariance and Einstein
causality on an arbitrary $T$ determines a covariant QSP up to unitary
equivalence uniquely if one requires the fulfillment for the family $(\pi
_{t})_{t\in T}$ of minimal conditions of normalization with respect to an
initial state space (state vector). The corresponding theorem, which
establishes the existence of the covariant QSP described axiomatically by
the causal correlation kernels and its uniqueness up to equivalence `almost
everywhere' (i.e., up to unitary equivalence of minimal modifications of the
representations of weakly equivalent processes), announced for arbitrary
observable algebra $\mathcal{B}$ in \cite{10}, plays the same role in the
operational statistical physics of open quantum systems as Wightman's
reconstruction theorem \cite{11} in quantum field theory or Kolmogorov's
fundamental theorem \cite{12} in the classical theory of random processes.

We also consider the question of determining a covariant QSP in the narrower
sense recently suggested by Accardi \cite{13} and Lewis \cite{14} and also
the relationship of the construction presented here with the noncausal
reconstruction of \cite{15} which uses all (not only chronologically
ordered) correlation functions. We show that these functions, which in
contrast to the causal functions form an infinite system even in the case of
a finite set $T$, can also be formally defined as causal functions on the
lexicographically ordered product $\mathbb{Z}\times T$. Although the process
indexed by the set $\mathbb{Z}\times T$ does not have a direct physical
interpretation in the real time $T$, its formal causal reconstruction leads
to a covariant family of representations $\pi _{i,t}=\pi _{t}$ independent
of $i\in \mathbb{Z}$, identical to the noncausal on $T$ representations $\pi
_{t}$ of the reconstruction of \cite{15}. Thus, the QSP reconstruction
theorem in the narrow sense \cite{15}, like Kolmogorov's classical
reconstruction \cite{12}, is a special case of the fundamental
reconstruction theorem formulated and proved here.

\section{Quantum Random Processes in the Narrow Sense}

1. Let $\mathcal{H}$ be a Hilbert space, $\mathcal{P}(\mathcal{H})$ denote
the set of Hermitian projection operators, $P^{\ast }=P=P^{2}$, on $\mathcal{%
H}$, $T$ be a linearly ordered set (e.g. the discrete or continuous time),
and with each $t\in T$ there be associated a Boolean $\sigma $-semiring $%
\mathcal{B}_{t}$. I.e. each $\mathcal{B}_{t}$ is a system of subsets $%
B\subseteq \mathrm{E}_{t}$ with the identity $\mathrm{E}_{t}\in \mathcal{B}%
_{t}$ and Boolean zero {$\mathrm{O}$}${_{t}}\in \mathcal{B}_{t}$ -- the
empty subset {$\mathrm{O}$}${_{t}}=\emptyset $ of $\mathrm{E}_{t}$; $%
\mathcal{B}_{t}$ is invariant with respect to the operation of
multiplication $BB^{\prime }=B\cap B^{\prime }$, and for any $B\in \mathcal{B%
}_{t}$ there exists a $\sigma $-partitioning $\left\{ B^{m}\right\}
_{1}^{\infty }\subseteq \mathcal{B}_{t}$ of the identity that contains $B$: $%
\mathrm{E}_{t}=\sum_{m=1}^{\infty }B^{m}$, $B\in \{B^{m}\}_{1}^{\infty }$ ($%
\sum_{m=1}^{\infty }B^{m}\;$denotes the union $\cup _{m=1}^{\infty }B^{m}$
of the disjoint events $B^{m}\in \mathcal{B}_{t}$, $B^{m}B^{m^{\prime }}=%
\mathrm{O}_{t}$ for$\;m\neq m^{\prime }$). As $\mathcal{B}_{t}$ for each $t$
one can take the smallest separating semiring of all single point subsets $%
B=\left\{ x\right\} $ with the total set $B=\mathrm{E}_{t}$ of a discrete
countable set $\mathrm{E}_{t}\subseteq \mathbb{N}$ or the largest one, the
power set Boolean algebra $\wp \left( \mathrm{E}_{t}\right) $, or the
semiring of intervals $B=[x,x^{\prime })$ of $\mathrm{E}_{t}\subseteq 
\mathbb{R}$, or the whole Borel $\sigma $-algebra $\mathcal{B}(\mathbb{R})$
generated by these intervals. The elements $B\in \mathcal{B}_{t}$ describe
the events `$x\in B$' at the time $t\in T$ with unit event $\mathrm{E}_{t}$
determining the readings $x\in \mathrm{E}_{t}$ of some measuring device with
the scale $\mathrm{E}_{t}$, which in general may depend on $t$. We shall
assume that on $T$ there acts a semigroup (or group) $S$ of symmetries such
that $st<st^{\prime }$ if $t<t^{\prime }$, this semigroup being represented
on $\mathcal{H}$ by linear isometries (or unitaries) $V_{s}:\mathcal{H}%
\rightarrow \mathcal{H}$, $V_{s}^{\ast }V_{s}=\mathrm{I}$, $s\in S$, and on
each $\mathrm{E}_{t}$ by measurable injections $g_{s}:\mathrm{E}%
_{t}\rightarrow \mathrm{E}_{st}$, $g_{s}g_{s^{\prime }}=g_{ss^{\prime }}$,
which determine Boolean homomorphisms $\mathcal{B}_{st}\ni B\mapsto
B^{s}=g_{s}^{-1}(B)$ onto $\mathcal{B}_{t}$ for every $t\in T$.

\begin{definition}
A quantum stochastic process over the family $\mathcal{B}=(\mathcal{B}_{t})$
in the narrow sense is described on the Hilbert space $\mathcal{H}$ by a
family $\pi =(\pi _{t})$ of $\sigma $-homomorphisms $\pi _{t}:\mathcal{B}%
_{t}\rightarrow {\mathcal{P}}${$\left( \mathcal{H}\right) $}, $t\in T$ and a
vector state specified by a normalized element $\xi \in \mathcal{H}$. The $%
\sigma $-homomorphisms are being defined as mappings of $B\in \mathcal{B}%
_{t} $ into the quantum logic of orthogonal projectors $P\in {\mathcal{P}}${$%
\left( \mathcal{H}\right) $} by the conditions 
\begin{equation}
\mathrm{E}_{t}=\sum_{m=1}^{\infty }B^{m},\quad B^{m}\in \mathcal{B}%
_{t}\Rightarrow \mathrm{I}=\sum_{m=1}^{\infty }P^{m},\quad P^{m}=\pi
_{t}(B^{m})
\end{equation}%
holding for each $t\in T$ ($B^{m}B^{m^{\prime }}=\mathrm{O}_{t}\Rightarrow
P^{m}P^{m^{\prime }}=0$ for $m\neq m^{\prime }$). The QSP is said to be $S$%
-stationary in the narrow sense if for every $s\in S$ and $t\in T$ the
conditions 
\begin{equation*}
V_{s}\pi _{t}(B^{s})=\pi _{st}(B)V_{s},\quad \forall B\in \mathcal{B}%
_{t},\quad t\in T
\end{equation*}%
and $V_{s}\xi =\xi $, $s\in S$, are satisfied with respect to a
representation $V=(V_{s})$ of $S$ on $\mathcal{H}$.
\end{definition}

Note that without loss of generality, every QSP can be considered as $S$-
stationary with respect to the given semigroup $S$ if one admits the trivial
action $t=st$ for all $s\in S$ on $T$ and $V_{s}=\mathrm{I}$ on $\mathcal{H}$%
.

As an example of $\sigma $-homomorphisms describing an $S$-stationary QSP
with respect to a given unitary representation $U=(U_{s})$ with invariant
state vector $\xi $, one can consider a family $\pi $ of the
projection-valued measures $\pi _{t}(B)=E_{t}(x^{\prime })-E_{t}(x)$ for $%
B=[x,x^{\prime })$, determined by spectral families $\{E_{t}(x),x\in \mathbb{%
R}\}$ of self-adjoint operators $X_{t}=\int xdE_{t}(x)$, $t\in T$ in $%
\mathcal{H}$ that transform covariantly with respect to a state vector $\xi $%
: 
\begin{equation*}
U_{s}X_{t}=X_{st}U_{s}\quad \forall t\in T,\;\;U_{s}\xi =\xi \quad \forall
s\in S.
\end{equation*}%
By such a family $X=(X_{t})$ of unitarily equivalent (in the case of
transitivity of $S$ on $T$) operators $X_{t}=X_{t}^{\ast }$ one can specify
any real QSP that is $S$-stationary with respect to the trivial
representation $g_{s}(x)=x$ of the group $S$ on $\mathrm{E}_{t}=\mathbb{R}$ $%
\forall t\in T$. However for time-like measurements, when $\mathrm{E}%
_{t}\subseteq T$, the stationarity condition must be determined with respect
to nontrivial transformations $g_{s}(x)=sx$ from $\mathrm{E}_{t}$ into $%
\mathrm{E}_{st}=s\mathrm{E}_{t}$ as it is in the case of the translations $%
st=t+s$ on the additive group $T=\mathbb{R}$. For example, in the case $%
\mathrm{E}_{t}=\{x<t\}$, corresponding to the measurement at each $t$ of the
occurrence times $x\in \mathbb{R}$ of almost surely past events like the
birth times of a historic phenomena or starting times of a continuous
measurements, a QSP translationally invariant with respect to an additive
semigroup $S\subseteq \mathbb{R}$ can be determined by the spectral
decompositions of bounded from above selfadjoint operators $X_{t}<t\mathrm{I}
$ that together with the vector $\xi \in \mathcal{H}$ satisfy the covariance
condition 
\begin{equation*}
X_{t}+s\mathrm{I}=U_{s}^{\ast }X_{t+s}U_{s}\quad \forall t\in T,\;\;U_{s}\xi
=\xi \quad \forall s\in S.
\end{equation*}

2. We denote by $\mathcal{F}$ the set of finite parts $\Lambda \subset T$,
with each $\Lambda =\{t_{1},\dots ,t_{n}\}\equiv \Lambda _{n}$ being a chain 
$t_{n}>\dots >t_{1}$ of length $n$, and we let $\mathcal{B}^{\Lambda
}=\times _{t\in \Lambda }\mathcal{B}_{t}$ for every $\Lambda =\Lambda _{n}$
be the set of sequences $\mathbf{b}=(B_{1},\dots ,B_{n})$ of events $\mathbf{%
b}(t_{i})=B_{i}\in \mathcal{B}_{t_{i}}$. Each such $\mathbf{b}\in \mathcal{B}%
^{\Lambda }$ is determined by a unique function $b:T\ni t\mapsto b(t)\in 
\mathcal{B}_{t}$, $b\left( t\right) =\mathrm{E}_{t}$, $t\notin \Lambda $ as
the restriction $\mathbf{b}=b|\Lambda $, and $\mathbf{e}=\left( \mathrm{E}%
_{1},\ldots ,\mathrm{E}_{n}\right) \in \mathcal{B}^{\Lambda }$ is determined
by the function $e(t)=\mathrm{E}_{t}$ for all $t\in T$. Note that the family 
$\{\mathcal{B}^{\Lambda }\}$ is inductive with respect to the extension $%
\mathbf{b}\in \mathcal{B}^{\Lambda }\mapsto \hat{b}_{\mathrm{M}}\in \mathcal{%
B}^{\mathrm{M}}$ such that $\hat{b}_{\mathrm{M}}(t)=\mathrm{E}_{t}$ for $%
t\in \mathrm{M}\backslash \Lambda $ on any $\mathrm{M}\supseteq \Lambda $, $%
\mathrm{M}\in \mathcal{F}$ , defined by the restriction $\hat{b}_{\mathrm{M}%
}=\hat{b}|\mathrm{M}$ of the corresponding function $\hat{b}$ on $T$ for $%
\mathbf{b}=\hat{b}|\Lambda $. The smallest element $\mathcal{B}^{\emptyset }$
of the inductive family $\{\mathcal{B}^{\Lambda }\}$ consists of a single
element --- the empty sequence which can be extended on any $\Lambda $ as $%
e_{\Lambda }$ by the identity function $e$ on $T$.

In accordance with the statistical interpretation of quantum mechanics, the
probability $\mu ^{\Lambda }(\mathbf{b})$ of successive observation of the
events $B=b(t)$, $t\in \Lambda $, in the corresponding chain $\Lambda _{n}$
is determined by the family $\pi $ and the vector $\xi $ in agreement with
the von Neumann projection postulate \cite{4} as 
\begin{equation}
\mu ^{\Lambda }(\mathbf{b})=\left\Vert P_{n}\cdots P_{1}\xi \right\Vert
^{2}=\Vert \xi ^{\Lambda }(\mathbf{b})\Vert ^{2}.
\end{equation}%
Here $\xi ^{\Lambda _{n}}(\mathbf{b})=\pi _{t_{n}}(B_{n})\cdots \pi
_{t_{1}}(B_{1})\xi $ is the result of the chronologically ordered action of
the orthogonal projectors $P_{i}=\pi _{t_{i}}(B_{i})$ on the state vector $%
\xi $. On the basis of these probabilities, which can be determined
experimentally by counting the relative frequencies of the occurrences of
the event sequences $\mathbf{b}=(B_{1},\dots ,B_{n})$ when they are measured
in a real flow of time $t\in \{t_{1},\dots ,t_{n}\}$ on each copy of the
quantum ensemble, one can calculate different characteristics of the QSP and
even attempt to reconstruct it by building a statistically equivalent
mathematical model $\left( \mathcal{H},\pi ,\xi \right) $ of the real
quantum system in which this process is observed.

Such QSP reconstruction problem, considered in the present paper, should
play a key testing role for the foundation of any quantum dynamical theory.
Its solution establishes necessary and sufficient conditions under which
there exists --- and is unique up to an equivalence --- a minimal
mathematical model of the physical system under consideration, the model
correctly predicting the statistics of successive measurements of
observations in this system as a quantum stochastic processes in the above
narrow or in a wider sense.

In the framework of classical theory, this problem was solved by
Kolmogorov's fundamental theorem \cite{12}, a necessary and sufficient
condition for the applicability of this theorem being that the family $\{\mu
^{\Lambda }\}$ must form a projective system of probability measures. It is
readily verified that the family of positive normalized mappings $\mu
^{\Lambda }:\mathcal{B}^{\Lambda }\rightarrow \lbrack 0,1]$ does indeed form
a projective system, $\mu ^{\Lambda }(\mathbf{b})=\mu ^{\mathrm{M}}(\mathbf{b%
})$ for $\Lambda \subseteq \mathrm{M}\in \mathcal{F}$, but these mappings
are not in general additive, although they satisfy the condition of $\sigma $%
-additivity with respect to the last argument $B=\mathbf{b}(t)$, $t=\max
\Lambda $. As with respect to the remaining arguments of $\mathbf{b}(t)$, $%
t<\max \Lambda $, the probabilities $\mu ^{\Lambda }$ are not in general
even finitely additive due to the possible noncommutativity of the
orthoprojectors $\{P_{i}\}$ corresponding to different $t_{i}\in \Lambda $
in accordance with the causal dependence of the chronologically ordered
events $\{B_{i}\}$. The absence of this additivity, observed experimentally
in quantum interference processes, indicates the inadequacy of the classical
probability theory dealing only with the additive measures corresponding to
the compatible events which can always be represented by the commuting
orthoprojectors, $P_{i}P_{k}=\pi \left( B_{i}B_{k}\right) =\pi \left(
B_{k}B_{i}\right) =P_{k}P_{i}$, even if $t_{i}\neq t_{k}$. Thus in the
framework of the classical theory it is impossible not only to construct an
adequate mathematical model of the physical system capable of predicting the
noncommutative QSP statistics, but even to formulate the problem of its
reconstruction due to the nonadequacy of quantum probabilities $\mu
^{\Lambda }$ to the (additive) probability measures.

3. For the QSP reconstruction, instead of the probabilities $\mu ^{\Lambda }$%
, it is necessary to use the correlation kernels (multikernels) 
\begin{equation}
\kappa ^{\Lambda }(\mathbf{b,b}^{\prime })=(\xi ^{\Lambda }(\mathbf{b})\mid
\xi ^{\Lambda }(\mathbf{b}^{\prime })),\quad \mathbf{b,\;b}^{\prime }\in 
\mathcal{B}^{\Lambda },\;\Lambda \in \mathcal{F}.
\end{equation}%
They are in principle determined by means of the polarization formulas from
the diagonal values $\kappa ^{\Lambda }(\mathbf{b,b})$ of these kernels,
extended to the multi-sesquilinear forms on all possible operator sequences $%
\mathbf{b}$. The multikernels $\kappa ^{\Lambda }:\mathcal{B}^{\Lambda
}\times \mathcal{B}^{\Lambda }\rightarrow \mathbb{C}$ corresponding to the
QSP in the narrow sense $\pi $ on $\mathcal{B}$ obviously form a projective
system: $\kappa ^{\Lambda }(\mathbf{b,b}^{\prime })=\kappa ^{\mathrm{M}}(%
\mathbf{b,b}^{\prime })$, if $\Lambda \subseteq \mathrm{M}\in \mathcal{F}$, $%
\mathbf{b,\;b}^{\prime }\in \mathcal{B}^{\Lambda }$, and for every $\Lambda
\in \mathcal{F}$ are readily described by the following properties 1--4,
which are called, respectively, positivity, normalisability, $\sigma $%
-additivity, and factorisability with respect to the product $(B,B^{\prime
})\mapsto BB^{\prime }$:

\begin{enumerate}
\item[1.] $\sum_{i,i^{\prime }=1}^{m}\kappa ^{\Lambda }(\mathbf{b}^{i},%
\mathbf{b}^{i^{\prime }})\bar{c}_{i}c_{i^{\prime }}\geq 0=$, $\forall 
\mathbf{b}^{i}\in \mathcal{B}^{\Lambda }$, $c_{i}\in \mathbb{C}$, $i\leq
m\!=\!1,2,\dots $.

\item[2.] $\kappa ^{\Lambda }(\mathbf{e,e})=1$, in particular $\kappa
^{\emptyset }=1$.

\item[3.] $\kappa ^{\Lambda }(\mathbf{b,b})=\sum_{m=1}^{\infty }\kappa
^{\Lambda }(\mathbf{b}^{m},\mathbf{b}^{m})$, where $\mathbf{b}^{m}(t)=%
\mathbf{b}(t)$, $t<t_{n}=\max \Lambda $, and 
\begin{equation*}
b(t_{n})=\sum_{m=1}^{\infty }b^{m}(t_{n})\quad (b^{m}(t_{n})b^{m^{\prime
}}(t_{n})=\mathrm{O}\;\;\text{for}\;\;m\neq m^{\prime }).
\end{equation*}

\item[4.] $\kappa ^{\Lambda }(\mathbf{b}B,\mathbf{b}^{\prime })=\kappa
^{\Lambda }(\mathbf{b},\mathbf{b}^{\prime }B)$ for any $B\in \mathcal{B}%
_{tn} $, $\mathbf{b},\;\mathbf{b}^{\prime }\in \mathcal{B}^{\Lambda }$,
where $(\mathbf{b}B)(t)=\mathbf{b}(t)$, $t\neq t_{n}$, $(\mathbf{b}%
B)(t_{n})=b(t_{n})B$, $t_{n}=\max \Lambda $.
\end{enumerate}

For a QSP that is $S$-stationary in the narrow sense, the multikernels $%
\kappa ^{\Lambda }$ also satisfy the condition of $S$-stationarity in the
wide sense:

\begin{enumerate}
\item[5.] $\kappa ^{s\Lambda }(\mathbf{b},\mathbf{b}^{\prime })=\kappa
^{\Lambda }(\mathbf{b}^{s},\mathbf{b}^{\prime s})$, $\forall \mathbf{b},%
\mathbf{b}^{\prime }\in \mathcal{B}^{s\Lambda }$, $s\in S$, where $s\Lambda
=\{st:t\in \Lambda \}$, $b^{s}(t)=b(st)^{s}$, $t\in \Lambda $.
\end{enumerate}

As follows from Theorem 2 (see Sec.2), the natural conditions 1--5 also hold
for the wider definition of QSPs described in Sec.2 by weakened conditions
of normalization $\pi _{t}(\mathrm{E}_{t})=E_{t}$, $E_{t}\xi =\xi $ instead
of the condition $\pi _{t}(\mathrm{E}_{t})=\mathrm{I}$, $\forall t\in T$,
these holding with respect to a nondecreasing family $E=(E_{t})_{t\in T}$ of
orthogonal projectors $E_{t}\leq E_{t^{\prime }}$, $t\leq t^{\prime }$. Such
widening of the QSP concept makes it possible to reconstruct the QSP from an
arbitrary projective system $\{\kappa ^{\Lambda }\}$ determined by Axioms
1--5 in a canonical way as it is formulated in the fundamental Theorem 3.
Theorem 4 shows that this widening of QSPs, which requires fulfillment of
the ordinary normalization condition only in the expanding system of
subspaces $\mathcal{H}_{t}=E_{t}\mathcal{H}\ni \xi $, is necessary if the
minimal process, defined in Sec.4, is to be described by the system $%
\{\kappa ^{\Lambda }\}$ uniquely (up to unitary equivalence).

4. We now show how the canonical QSP reconstruction presented for an
arbitrary preodered set $\mathrm{T}$ in Sec.3 makes it possible to obtain
the canonical reconstruction of a QSP in the narrow sense as the particular
case of the Theorem 3, using the much larger system of \emph{all}, not only
causal multikernels $\kappa _{\boldsymbol{x}}(\mathbf{b},\mathbf{b}^{\prime
})$, indexed by arbitrary finite sequences $\boldsymbol{x}\in \cup
_{n=0}^{\infty }X^{n}$. Such reconstruction, recently suggested in \cite{15}%
, requires the considering of not necessarily chronological orderings $\pi
_{x_{1}}(B_{1})\cdots \pi _{x_{n}}(B_{n})\xi =\xi _{\boldsymbol{x}}(\mathbf{b%
})$ of the events $\mathbf{b}=(B_{1},\dots ,B_{n})$ since the elements of a
sequence $\boldsymbol{x}=(x_{1},\dots ,x_{n})\in X^{n}$ may not form a chain 
$x_{1}>\dots >x_{n}$. Therefore such kernels cannot be not physically
realized by the sequential measurements in real time unless points $x\in X$
are considered as causally equivalent, or merely disordered (space-like)
coordinates $x\in \mathbb{R}^{d}$, in which case the families $\pi =\left(
\pi _{x}\right) $ give kinematic representation of quantum fields but not
dynamic processes. The following construction should be therefore considered
as a reconstruction of quantum fields rather than the processes.

Let $\mathrm{T}=\mathbb{Z}\times X$ be the set of pairs $\mathrm{t}=(l,x)$, $%
l=0,\pm 1,\pm 2,\dots $, linearly preordered as $\mathrm{t}\lesssim \mathrm{t%
}^{\prime }\Leftrightarrow $ $l\leq l^{\prime }$. It can be completely
ordered by the lexicographic order $\mathrm{t}\leq \mathrm{t}^{\prime
}\Leftrightarrow $ $x\leq x^{\prime }$ if $l=l^{\prime }$, or we can leave
such points causally equivalent, $\mathrm{t}\sim \mathrm{t}^{\prime }$ if $%
l=l^{\prime }$, and identify with the corresponding points $x\sim x^{\prime
}\in X$ for $l=0$. Let $S=\mathbb{Z}$ be the group of translations $\mathrm{t%
}=(l,x)\mapsto s\mathrm{t}=(l+s,x)$, acting on $\mathrm{T}$
quasitransitively in the following sense: For any pair $\mathrm{t}>\mathrm{t}%
^{\prime }$, there exists $s\in \mathbb{Z}$ such that $\mathrm{t}^{\prime }>s%
\mathrm{t}$ (it is necessary to take $s<l^{\prime }-l$).

We consider the QSP $\pi _{l,x}=\pi _{x}$ over $\mathcal{B}_{l,x}=\mathcal{B}%
_{x}$, $\left( l,x\right) =\mathrm{t}\in \mathrm{T}$, determined by a QSP
process $\pi _{x}:\mathcal{B}_{x}\rightarrow \mathcal{P}\left( \mathcal{H}%
\right) $ on $X\ni x$ in the narrow sense on the Hilbert space $\mathcal{H}$
with state vector $\xi \in \mathcal{H}$. Such defined process $\pi _{\mathrm{%
t}},\mathrm{t}\in \mathrm{T}$ is $S$-stationary in the narrow sense due to
the condition of quasiconstancy 
\begin{equation}
\pi _{\mathrm{t}}(B)=\pi _{s\mathrm{t}}(B),\quad \forall B\in \mathcal{B}_{%
\mathrm{t}},\quad \mathrm{t}\in \mathrm{T},
\end{equation}%
which is the $S$-covariance with respect to the trivial representations $%
g_{s}=\mathrm{id}$, $U_{s}=1$ of the group $S$ on $\mathrm{E}_{\mathrm{t}}=%
\mathrm{E}_{x}$ and $\mathcal{H}$. It is readily verified that this
quasiconstant process induces a system of multikernels (1.3) indexed by
pairs $\Lambda =\left( \boldsymbol{s},\boldsymbol{x}\right) $ of $%
\boldsymbol{s}=(s_{1},\dots ,s_{n})\in \mathbb{Z}^{n}$, $s_{1}\leq \ldots
\leq s_{n}$, and $\boldsymbol{x}\in X^{n}$, satisfying the strengthened
condition of $S$-stationarity (ultrastationarity):

\begin{enumerate}
\item[5'.] $\kappa _{\boldsymbol{x}}^{\boldsymbol{s}}(\mathbf{b},\mathbf{b}%
^{\prime })=\kappa _{\boldsymbol{x}}^{\boldsymbol{0}}(\mathbf{b},\mathbf{b}%
^{\prime })$\ $\forall \mathbf{b,b}^{\prime }\in \mathcal{B}^{\Lambda
}=\times _{i=1}^{n}\mathcal{B}_{x_{i}}$, where $\boldsymbol{0}=\left(
0,\ldots ,0\right) $.
\end{enumerate}

By virtue of condition $5^{\prime }$, the projective system $\{\kappa
^{\Lambda }=\kappa _{\boldsymbol{x}}^{\boldsymbol{s}}\}$, indexed by finite $%
\Lambda \subset \mathrm{T}$, is determined by a set $\{\kappa _{\boldsymbol{x%
}}^{\mathbf{0}}=\kappa _{\boldsymbol{x}}\}$ of arbitrarily ordered
multikernels $\kappa _{\boldsymbol{x}}$ of the original process (field) $\pi
_{x}$, $x\in X$, these being identical for any $\boldsymbol{x}=(x_{1},\dots
,x_{n})$ to $\kappa ^{\Lambda }$ with $s_{1}<\dots <s_{n}$ and thus with $%
\Lambda =\left\{ \mathrm{t}_{1}<\dots <\mathrm{t}_{n}\right\} $.

\begin{theorem}
Let $\{\kappa ^{\Lambda }\}$ be a projective system of multikernels $\kappa
^{\Lambda }:\;\mathcal{B}^{\Lambda }\times \mathcal{B}^{\Lambda }\rightarrow 
\mathbb{C}$ that satisfy conditions 1--4 for every finite chain $\Lambda
\subset \mathrm{T}$ and also condition $5^{\prime }$ of ultrastationarity
with respect to the quasitransitive action of the group $S=\mathbb{Z}$ on $%
\mathrm{T}$. Then there exists a Hilbert space $\mathcal{H}$, a family of $%
\sigma $-homomorphisms $\pi _{\mathrm{t}}:\mathcal{B}_{\mathrm{t}%
}\rightarrow \mathcal{P}\left( \mathcal{H}\right) $, $\pi _{\mathrm{t}}(%
\mathrm{E}_{\mathrm{t}})=\mathrm{I}$, $\mathrm{t}\in \mathrm{T}$, and a
state vector $\xi \in \mathcal{H}$ that determine an $S$-quasiconstant QSP $%
\pi =(\pi _{\mathrm{t}})$ in the narrow sense over $\mathcal{B}=(\mathcal{B}%
_{\mathrm{t}})$ which induces the system $\{\kappa ^{\Lambda }\}$. On the
minimal space $\mathcal{H}^{\ast }$ determined by the condition of
separability of the commutant $\mathcal{A}=\cap \pi _{l}(\mathcal{B}%
)^{\prime }$ with respect to $\xi $, this process is described uniquely up
to unitary equivalence.
\end{theorem}

The\/ proof of this theorem, which is a corollary to Theorem 3, is obtained
in (see Sec.\ 3) by the construction of the canonical process $\pi _{\mathrm{%
t}}^{\ast }$, $\mathrm{t}\in \mathrm{T}$, on the Hilbert space $\mathcal{H}$
that determines the minimal decomposition (1.3). By virtue of $5^{\prime }$,
this process is quasiconstant with respect to the quasitransitive action of $%
S$ on $\mathrm{T}$ and has constant unit $E_{l}=\pi _{l}^{\ast }(\mathrm{E}%
_{l})=E_{0}$ for any $l\in \mathbb{Z}$ due to $E_{l}\leq E_{l^{\prime }}\leq
E_{l}$ for any pair $l<l^{\prime }$ of the ordered set $\mathbb{Z}$ as will
be shown in Sec.\ 3. by virtue of the monotonicity of the family $E_{l}$, $%
l\in \mathbb{Z}$. Therefore, $E_{l}=\mathrm{I}$ on the minimal space $%
\mathcal{H}$ on which this process is uniquely determined by virtue of
Theorem 4. For the case $\mathrm{T}=\mathbb{Z}\times X$, the described
canonical reconstruction determines from the family $\{\kappa _{\boldsymbol{x%
}}\}=\{\kappa _{\boldsymbol{x}}^{\mathbf{0}}\}$ the QSP $\pi _{x}=\pi _{%
\mathrm{t}}$, $\mathrm{t}=\left( l,x\right) $, in the narrow sense.

\section{Definition and Properties of a General QSP}

1. Let $\mathrm{T}$ be an arbitrary set with elements $\mathrm{t}\in \mathrm{%
T}$ that is preordered by the reflexive and transitive relation $\leq $,
which is called causality on $\mathrm{T}$, and let $S$ be a semigroup (or
group) of causality preserving symmetries $s:\mathrm{T}\rightarrow \mathrm{T}
$ as monotonic transformations $\mathrm{t}\leq \mathrm{t}^{\prime
}\Leftrightarrow s\mathrm{t}\leq s\mathrm{t}^{\prime }$ (and therefore $%
\mathrm{t}\sim \mathrm{t}^{\prime }$ if $s\mathrm{t}=s\mathrm{t}^{\prime }$%
). We denote by the symbol $\sim $ the equivalence relation $\mathrm{t}\sim 
\mathrm{t}^{\prime }\Leftrightarrow \mathrm{t}\leq \mathrm{t}^{\prime }\leq 
\mathrm{t}$, which means invertibility of the causal dependence of $\mathrm{t%
}$ and $\mathrm{t}^{\prime }\in \mathrm{T}$, by $\Join $ the
noncomparability relation $\mathrm{t}\Join \mathrm{t}^{\prime
}\Leftrightarrow \mathrm{t}\not\leq \mathrm{t}^{\prime }\not\leq \mathrm{t}$%
, which determines causal independence of the corresponding $\mathrm{t}$ and 
$\mathrm{t}^{\prime }$, and by $>$ the strict order (anticipation) relation $%
\mathrm{t}>\mathrm{t}^{\prime }\Leftrightarrow \mathrm{t}\not\sim \mathrm{t}%
^{\prime }\leq \mathrm{t}$, and we shall say that a pair $\mathrm{t}$, $%
\mathrm{t}^{\prime }$ is mutually nonanticipatory if $\mathrm{t}\ngtr 
\mathrm{t}^{\prime }\ngtr \mathrm{t}$, i.e., if $\mathrm{t}\sim \mathrm{t}%
^{\prime }$, or $\mathrm{t}\Join \mathrm{t}^{\prime }$.

Besides linearly ordered subsets $\mathrm{T}\subseteq \mathbb{R}$, it may be
of interest to take as such $\mathrm{T}$ any space-time region of $(d+1)$%
-dimensional Minkowski space of $\mathrm{t}=(\tau ,\mathbf{r})$, where $\tau
\in \mathbb{R}$, $\mathbf{r}\in \mathbb{R}^{d}$, partially ordered by the
relation $\mathrm{t}\leq \mathrm{t}^{\prime }\Leftrightarrow c^{-1}|\mathbf{r%
}-\mathbf{r}^{\prime }|\leq \tau ^{\prime }-\tau $ which has trivial
equivalence $\mathrm{t}\sim \mathrm{t}^{\prime }\Leftrightarrow \mathrm{t}=%
\mathrm{t}^{\prime }$ and nonempty (for $d\neq 0$) Einsteinian causal
independence $\mathrm{t}\Join \mathrm{t}^{\prime }\Leftrightarrow |\mathbf{r}%
-\mathbf{r}^{\prime }|>c|\tau -\tau ^{\prime }|$. As $S$ one can consider
any semigroup (subgroup) of the inhomogeneous Lorentz transformations $s%
\mathrm{t}=\Lambda \mathrm{t}+\mathrm{l}$ that preserves the order and leave
the subset $\mathrm{T}\subseteq \mathbb{R}^{d+1}$ invariant. For example, if 
$\mathrm{T}$ is the future cone $\mathrm{t}\geq 0:\;c\tau \geq |\mathbf{r}|$%
, then as $S$ we can consider the semigroup of translations $\mathrm{t}%
\mapsto \mathrm{t}+\mathrm{l}$, $\mathrm{l}\in \mathrm{T}$, and the proper
orthochronous Lorentz group $L_{+}^{\uparrow }$. Galilean causality,
corresponding to the limit $c\rightarrow \infty $, is described, in
contrast, by the preorder relation $\mathrm{t}\leq \mathrm{t}^{\prime
}\Leftrightarrow \tau \leq \tau ^{\prime }$ with nontrivial (for $d\neq 0$)
equivalence $\mathrm{t}\sim \mathrm{t}^{\prime }\Leftrightarrow \tau =\tau
^{\prime }$ and empty noncomparability relation $\Join $, i.e., it is a
linear preorder given by the foliation of $\mathbb{R}^{d+1}$ into the
hyperplanes $\tau =\mathrm{const}$. At the same time, any spatial region $%
\mathrm{T}$ corresponding to fixed $\tau \in \mathbb{R}$ is equipped with
trivial causality: $\mathrm{t}\sim \mathrm{t}^{\prime }$ for any pair $%
\mathrm{t}$, $\mathrm{t}^{\prime }\in \mathrm{T}$, in contrast to the
relativistic case $c<\infty $, for which any spatial region has identical
causality $\mathrm{t}\leq \mathrm{t}^{\prime }\Rightarrow \mathrm{t}=\mathrm{%
t}^{\prime }$.

We introduce the following notation: $J$ is the class of all subsets $%
j\subset \mathrm{T}$ of pairwise \emph{nonanticipatory} elements $\mathrm{t},%
\mathrm{t}^{\prime }\in j\Rightarrow \mathrm{t}\ngtr \mathrm{t}^{\prime
}\ngtr \mathrm{t}$; $K$ is the subclass formed by the \emph{finite} subsets $%
k\in J$; $L$ is the subclass formed by the maximal subsets $l\in J$ such
that $l\subseteq j\in J\Rightarrow l=j$; and $T\subset J$ is the factor-set
of $\mathrm{T}$ formed by the maximal subsets $t\subset \mathrm{T}$ of
pairwise \emph{equivalent} elements $\mathrm{t}$, $\mathrm{t}^{\prime }\in
t\Rightarrow \mathrm{t}\sim \mathrm{t}^{\prime }$. The class $J$, whose
elements can be determined by the condition $j=\max j$, where $\max j$ is
the subset of the elements $\mathrm{t}\in j$ which are maximal in $j$ in the
sense $\mathrm{t}\leq \mathrm{t}^{\prime }\in j\Rightarrow \mathrm{t}\sim 
\mathrm{t}^{\prime }$, contains together with $j\in J$ any subset of it, $%
j^{\prime }\subseteq j$, and is a semilattice with respect to the operation $%
j\vee j^{\prime }=\max j\cup j^{\prime }$, which determines a strict order $%
j>j^{\prime }\Leftrightarrow j\vee j^{\prime }=j$, $j\cap j^{\prime
}=\emptyset $ of the partial preorder: 
\begin{equation*}
j\leq j^{\prime }\Leftrightarrow \;\forall \mathrm{t}\in j\;\;\exists 
\mathrm{t}^{\prime }\in j^{\prime }:\mathrm{t}\leq \mathrm{t}^{\prime }.
\end{equation*}%
Both these properties are also inherited by the subclass $K$, whereas the
subclass $L$ is only a directed set: any pair $l,l^{\prime }\in L$ has in $L$
a majorant $l_{+}\supseteq \max l\cup l^{\prime }$ and a minorant $%
l_{-}\supseteq \min l\cup l^{\prime }$. (It is assumed that for any $j\in J$
there exists $l\in L$ such that $j\subseteq l$). The factor-set $T$ is in
general an arbitrary partially ordered set. Subsets $j,j^{\prime }\in J$ are
said to be \emph{equivalent}: $j\sim j^{\prime }$, if $j\leq j^{\prime }\leq
j$; \emph{independent}: $j\Join j^{\prime }$, if $\mathrm{t}\Join \mathrm{t}%
^{\prime }$ for all $\mathrm{t}\in j$, $\mathrm{t}^{\prime }\in j^{\prime }$%
; and \emph{mutually nonanticipatory} if $j\vee j^{\prime }=j\cup j^{\prime
} $. In these terms the empty subset $\emptyset \subset \mathrm{T}$ is the
least element $\emptyset \in J:\;\emptyset \leq j\in J$; $j\sim \emptyset
\Rightarrow j=\emptyset $; $j\Join \emptyset $ for any $j\in J$; and $%
j>\emptyset $, if $j\neq \emptyset $.

2. Let $\mathfrak{A}=(\mathfrak{A}_{k})_{k\in K}$ be a nonincreasing family
of $\ast $-subalgebras of the C*-algebra $\mathfrak{A}_{\emptyset }=\mathcal{%
B}\left( \mathcal{K}\right) $ of bounded operators $a:\mathcal{K}\rightarrow 
\mathcal{K}$ of the Hilbert space $\mathcal{K}$ with common unit $1_{k}=1\in 
\mathfrak{A}_{k}$ such that $\mathfrak{A}_{k\vee k^{\prime }}=\mathfrak{A}%
_{k}\cap \mathfrak{A}_{k^{\prime }}$ for all $k,k^{\prime }\in K$, and $%
\mathcal{B}=(\mathcal{B}_{k})_{k\in K}$ be the family of Boolean semirings
of subsets $B\subseteq \mathrm{E}_{k}$ with units $\mathrm{E}_{k}$, $k\in K$%
, satisfying the conditions 
\begin{eqnarray}
k &\sim &k^{\prime }\Rightarrow \mathcal{B}_{k\cup k^{\prime }}\supseteq 
\mathcal{B}_{k}\cup \mathcal{B}_{k^{\prime }},\qquad \mathrm{E}_{k}=\mathrm{E%
}_{k^{\prime }}, \\
k &\Join &k^{\prime }\Rightarrow \mathcal{B}_{k\cup k^{\prime }}=\mathcal{B}%
_{k}\otimes \mathcal{B}_{k^{\prime }},\qquad \mathrm{E}_{k\cup k^{\prime }}=%
\mathrm{E}_{k}\times \mathrm{E}_{k^{\prime }},
\end{eqnarray}%
where $\mathcal{B}_{k}\otimes \mathcal{B}_{k^{\prime }}$ is a semiring of
subsets $B\times B^{\prime }$, $B\in \mathcal{B}_{k}$, $B^{\prime }\in 
\mathcal{B}_{k^{\prime }}$, of the Cartesian product $\mathrm{E}_{k}\times 
\mathrm{E}_{k^{\prime }}$ ($\mathcal{B}_{\emptyset }$ is identified in
accordance with (2.2) with the trivial Boolean algebra $(\mathrm{O},\mathrm{E%
})$ of a single-point set $\mathrm{E}_{\emptyset }$). Denoting $\mathfrak{A}%
_{t}=\mathfrak{\mathfrak{A}}_{\{\mathrm{t}\}}$ for any single-point subset $%
\{\mathrm{t}\}\subseteq t\in T$ and $\mathcal{B}_{t}=\underset{k\subset t}{%
\cup }\mathcal{B}_{k}$, $\mathrm{E}_{t}=\mathrm{E}_{\{\mathrm{t}\}}$, we
have for $k\neq \emptyset $: $\mathfrak{\mathfrak{A}}_{k}\subseteq \underset{%
t\in \lbrack k]}{\cap }\mathfrak{\mathfrak{A}}_{t}$, and $B=\underset{t\in
\lbrack k]}{\times }b(t)$ for any $B\in \mathcal{B}_{k}$, where $\mathrm{E}%
_{t}\supseteq b(t)\in \mathcal{B}_{t}$, $[k]=\{t\in T:t\cap k\neq \emptyset
\}$.

We assume that the semigroup $S$, acting on $K$ by transformations $sk=\{s%
\mathrm{t}:\mathrm{t}\in k\}$, is represented on $\mathfrak{A}$ by
self-consistent family of $\ast $-endomorphisms $a\mapsto a^{s}=U_{s}^{\ast
}aU_{s}$, where $U=(U_{s})_{s\in S}$ is a unitary representation of $S$ on $%
\mathcal{K}$, the endomorphisms mapping each subalgebra $\mathfrak{A}_{sk}$
onto $\mathfrak{A}_{k}$, and we assume that the semigroup is also
represented on $\mathcal{B}$ by $\sigma $-homomorphisms from each semiring $%
\mathcal{B}_{sk}$ onto $\mathcal{B}_{k}$, these being induced by measurable
mappings $g_{s}:\mathrm{E}_{t}\rightarrow \mathrm{E}_{st}$.

Let $\mathcal{H}$ be a Hilbert space containing $\mathcal{K}$, and $\iota
=(\iota _{k})_{k\in \mathbb{Z}}$ be a family of faithful $\ast $%
-representations $\iota _{k}:\mathfrak{A}_{k}\rightarrow \mathcal{B}\left( 
\mathcal{H}\right) $ that satisfy the self-consistency condition 
\begin{equation}
\iota _{k}(a)=\iota _{k^{\prime }}(a)I_{k},\quad \forall a\in \mathfrak{A}%
_{k^{\prime }},\;k\leq k^{\prime }\in K,
\end{equation}%
where $I_{k}=\iota _{k}(1)$ are orthoprojectors in $\mathcal{H}$, $\iota
_{\emptyset }$ is the identity representation on the subspace $\mathcal{K}%
=I_{\emptyset }\mathcal{H}$. We consider also an isometric representation $%
V=\left( V_{s}\right) _{s\in S}$ of the semigroup $S$ on $\mathcal{H}$ with
respect to which the condition of covariance of the representations $\iota $
is satisfied in the form 
\begin{equation}
V_{s}\iota _{k}(a^{s})=\iota _{sk}(a)V_{s}I_{k}\;\;\quad \forall a\in 
\mathfrak{A}_{sk}\;,
\end{equation}%
so that $U_{s}=V_{s}I_{\emptyset }$ for any $s\in S$. Note that the family $%
I=(I_{k})_{k\in K}$ of orthogonal projectors $I_{k}$ determining essential
subspaces $\mathcal{H}_{k}=I_{k}\mathcal{H}$ of the representations $\iota
_{k}$, is nondecreasing: $I_{k}=I_{k}I_{k^{\prime }}\leq I_{k^{\prime }}$
for any $k\leq k^{\prime }$, and uniquely determines these representations
in the case of a one-dimensional $\mathcal{K}\simeq \mathbb{C}$ \ (and
therefore $\mathfrak{A}_{k}=\mathbb{C}$ for all $k\in K$): $\iota
_{k}(c)=cI_{k}$, $c\in \mathbb{C}$ (at the same time $I_{\emptyset }$ is the
one-dimensional orthogonal projector $\mathrm{P}_{\xi }$ corresponding to a
state vector $\xi \in \underset{k>\emptyset }{\cap }\mathcal{H}_{k}$).

\begin{definition}
An $\mathcal{H}$-process with respect to $(\mathfrak{A},\iota )$ over $%
\mathcal{B}$ (or $\mathcal{H}$-QSP with respect to $I=(I_{k})$ if $\mathcal{K%
}\simeq \mathbb{C}$) is a family $\pi =(\pi _{k})_{k\in K}$ representations $%
\pi _{k}:\mathcal{B}\rightarrow \mathcal{P}\left( \mathcal{H}\right) $ that
satisfy the normalization conditions $\pi _{k}(E_{k})I_{k}=I_{k}$, $\forall
k\in K$, or even stronger conditions

\begin{enumerate}
\item[0)] $\underset{k\subseteq l}{\wedge }P_{k}=\underset{k\subseteq l}{%
\vee }I_{k}$, $\forall l\in L$, where%
\begin{equation*}
P_{k}=\pi _{k}(\mathrm{E}_{k}),\;\;I_{k}=\iota _{k}(1)\;\;\;\;\;\ \forall
k\in K,
\end{equation*}
where $P_{\emptyset }=\mathrm{I}$ is the identity in $\mathcal{H}$, $%
I_{\emptyset }=\mathrm{P}$ is the projection onto $\mathcal{K}$. These
representations are described by the following axioms:

\item[1)] $k\sim k^{\prime }\Rightarrow \pi _{k\cup k^{\prime }}(BB^{\prime
})=\pi _{k}(B)\pi _{k^{\prime }}(B^{\prime })$, $\forall B\in \mathcal{B}%
_{k} $, $B^{\prime }\in \mathcal{B}_{k^{\prime }}$,

\item[2)] $k\Join k^{\prime }\Rightarrow \pi _{k\cup k^{\prime }}(B\times
B^{\prime })=\pi _{k}(B)\pi _{k^{\prime }}(B^{\prime })P_{k\cup k^{\prime }}$%
,

\item[3)] If $B^{m}B^{m^{\prime }}=\mathrm{O}_{k}$ for $B^{m}\in \mathcal{B}%
_{k}$ with$\;m\neq m^{\prime }$, 
\begin{equation*}
\mathrm{E}_{k}=\sum_{m=1}^{\infty }B^{m}\Rightarrow P_{k}=\sum_{m=1}^{\infty
}\pi _{k}(B^{m}),
\end{equation*}

\item[4)] $[\pi _{k}(B),\iota _{k}(a)]=0$, $\forall B\in \mathcal{B}_{k}$, $%
a\in \mathfrak{A}_{k}$, $k\in K$.

The process $(\mathcal{B},\pi )$ is said to be $S$-covariant with respect to
the representation $V$ if for any $s\in S$

\item[5)] $V_{s}\pi _{k}(B^{s})=\pi _{sk}(B)V_{s}P_{k}$, $\forall B\in 
\mathcal{B}_{sk}$, $k\in K$.
\end{enumerate}
\end{definition}

3. By virtue of the uniqueness of the partitioning $\Lambda =\cup k_{i}$ of
any finite subset $\Lambda \subset T$ into elements $k_{i}\in K$ that form a
chain $k_{n}>\dots >k_{1}>\emptyset $ of a length $n\leq |\Lambda |$, we can
uniquely associate with every function $[\Lambda ]\ni t\mapsto $ $b(t)\in 
\mathcal{B}_{\Lambda \cap t}$, defined on a subset $[\Lambda ]\subset T$ of
equivalence classes $t=\left[ \mathrm{t}\right] $ such that $t\cap \Lambda
\neq \emptyset $, a sequence $\mathbf{b}=(B_{1},\dots ,B_{n})$ of events $%
B_{i}=\prod_{t\in \lbrack k_{i}]}b(t)=\mathbf{b}(k_{i})$. Here $k_{i}=\max
\Lambda _{i}$, $\Lambda _{i}\subseteq \Lambda $ are determined recursively: $%
\Lambda _{i-1}=\Lambda _{i}\backslash k_{i}$, $\Lambda _{n}=\Lambda $, and $%
n=n\left( \Lambda \right) $ is such that $\Lambda _{1}=k_{1}$, $\Lambda
_{0}=\emptyset $. For any $\Lambda \in \mathcal{F}$ such that $\Lambda >k$
for a $k\in K$, we denote by $\mathcal{B}_{k}^{\Lambda }$ the set of
sequences $\mathbf{b}$ as functions $B_{i}=\mathbf{b}(k_{i})$ on the chain $%
\left\{ k_{n}>\dots >k_{1}\right\} >k$ corresponding to such $\Lambda $, and
we determine an induction mapping $\mathbf{b}\in \mathcal{B}_{j}^{\Lambda
}\mapsto \hat{\mathbf{b}}\in \mathcal{B}_{k}^{\mathrm{M}}$ for any $j\in K$, 
$k\leq j$, $\Lambda \subseteq \mathrm{M}\in \mathcal{F}$, by extending the
function $b(t)$ corresponding to the sequence $\mathbf{b}$ to the function
on $\mathrm{M}$ as $\hat{b}(t)=b(t)$ for $t\in \lbrack \Lambda ]$ and $\hat{b%
}(t)=\mathrm{E}_{t}$ for $t\in \lbrack \mathrm{M}\backslash \Lambda ]$. This
function determines the sequence $\hat{\mathbf{b}}:\hat{\mathbf{b}}(k_{i})=%
\underset{t\in \lbrack k_{i}]}{\times }\hat{b}(t)$ with $i=1,\ldots ,n$, $%
n=n(\mathrm{M})$ such that $\hat{\mathbf{b}}(k_{i})=\mathbf{b}(k_{i})\times 
\mathrm{E}_{\bar{\Lambda}k_{i}}$, where $\mathbf{b}(k_{i})=\underset{t\in
\lbrack \Lambda k_{i}]}{\times }b(t)$, $\Lambda k_{i}=\Lambda \cap k_{i}$.

With each pair $\mathbf{b},\mathbf{b}^{\prime }\in \mathcal{B}_{k}^{\Lambda
} $ we associate the correlation operator 
\begin{equation}
\kappa _{k}^{\Lambda }(\mathbf{b},\mathbf{b}^{\prime })=F_{k}^{\Lambda }(%
\mathbf{b})^{\ast }F_{k}^{\Lambda }(\mathbf{b}^{\prime }),
\end{equation}%
where 
\begin{equation}
F_{k}^{\Lambda _{n}}(\mathbf{b})=\pi _{k_{n}}(B_{n})\cdots \pi
_{k_{1}}(B_{1})I_{k}
\end{equation}%
is the `Feynman integral" over the paths in $\left\{ b\left( t\right)
\right\} $ determined by the $\mathcal{H}$-process $\pi $ over $\mathcal{B}$%
. The following theorem establishes the properties of the operator-valued
multikernels $\kappa _{k}^{\Lambda }$, these properties when $k=\emptyset $, 
$\mathcal{K}\simeq \mathbb{C}$ being identical to the corresponding
properties of the scalar kernels (1.3) of the QSP with respect to $\mathfrak{%
A}_{t}=\mathbb{C}$, $t\in T$, defined in the narrow sense.

\begin{theorem}
Let $(\mathcal{B},\pi )$ be an $\mathcal{H}$-process with respect to $(%
\mathfrak{A},\iota )$ that is $S$-covariant with respect to $V$. Then
Eq.(2.5) determines mappings $\kappa _{k}^{\Lambda }:\mathcal{B}%
_{k}^{\Lambda }\times \mathcal{B}_{k}^{\Lambda }\rightarrow \iota _{k}(%
\mathfrak{A}_{l})^{\prime }$, where $\Lambda >k$, $l=\max \Lambda $, that
satisfy the following conditions:

\begin{enumerate}
\item[0)] (self-consistency) for any $\mathbf{b},\mathbf{b}^{\prime }\in 
\mathcal{B}_{j}^{\Lambda }$ 
\begin{equation*}
I_{k}\kappa _{j}^{\Lambda }(\mathbf{b},\mathbf{b}^{\prime })I_{k}=\kappa
_{k}^{\mathrm{M}}(\hat{\mathbf{b}},\hat{\mathbf{b}}^{\prime }),\;\;k\leq
j\in K,\;\Lambda \subseteq \mathrm{M}\in \mathcal{F};
\end{equation*}

\item[1)] (positivity) for any $\mathbf{b}^{i}\in \mathcal{B}_{k}^{\Lambda }$%
, $\eta _{i}\in \mathcal{H}_{k}$, $i\leq m$ 
\begin{equation*}
\sum_{i,i^{\prime }=1}^{m}(\eta _{i}\mid \kappa _{k}^{\Lambda }(\mathbf{b}%
^{i},\mathbf{b}^{i^{\prime }})\eta _{i^{\prime }})\geq 0,\quad m=1,2,\dots
;\Lambda >k;
\end{equation*}

\item[2)] (normalizability) $\kappa _{k}^{\Lambda }(\mathbf{e},\mathbf{e}%
)=I_{k}$, $\Lambda >k\in K$;

\item[3)] ($\sigma $-additivity [In the weak operator topology and also in
the sense of order convergence.]) 
\begin{equation*}
\mathbf{b}(l)=\sum_{m=1}^{\infty }\mathbf{b}^{m}(l),\quad l=\max \Lambda
\Rightarrow \kappa _{k}^{\Lambda }(\mathbf{b},\mathbf{b})=\sum_{m=1}^{\infty
}\kappa _{k}^{\Lambda }(\mathbf{b}^{m},\mathbf{b}^{m}),
\end{equation*}%
where $\mathbf{b}^{m}(k_{i})=\mathbf{b}(k_{i})$, $k_{i}\subseteq \Lambda
\backslash l$;

\item[4)] (factorizability) $\kappa _{k}^{\Lambda }(\mathbf{b}B,\mathbf{b}%
^{\prime })=\kappa _{k}^{\Lambda }(\mathbf{b},\mathbf{b}^{\prime }B)\!\!,$
where%
\begin{equation*}
(\mathbf{b}B)(l)=\mathbf{b}(l)B\!\!,(\mathbf{b}B)(k_{i})=\mathbf{b}%
(k_{i}),k_{i}\subseteq \Lambda \backslash l\;;
\end{equation*}

\item[5)] (covariance) for any $s\in S$, $\mathbf{b},\mathbf{b}^{\prime }\in 
\mathcal{B}_{k}^{\Lambda }$ 
\begin{equation*}
V_{k,s}^{\ast }\kappa _{sk}^{s\Lambda }(\mathbf{b},\mathbf{b}^{\prime
})V_{k,s}=\kappa _{k}^{\Lambda }(\mathbf{b}^{s},\mathbf{b}^{\prime s}),
\end{equation*}%
where $s\Lambda =\{st:t\in \Lambda \}$, $\mathbf{b}^{s}(l)=\mathbf{b}%
(sl)^{s} $, $\mathbf{b}^{\prime s}(l)=\mathbf{b}^{\prime }(sl)^{s}$, and $%
V_{k,s}:\mathcal{H}_{k}\rightarrow \mathcal{H}_{sk}$ are isometric
operators, made consistent by the condition $V_{k,s}|\mathcal{H}_{k^{\prime
}}=V_{k^{\prime },s}$, $k^{\prime }\leq k$, which form a representation, $%
V_{sk,s^{\prime }}V_{k,s}=V_{k,s^{\prime }s}$, of the semigroup $S$ on the
family $\mathcal{H}_{k}=I_{k}\mathcal{H}$, $k\in K$, with unitary
subrepresentation $U_{s}=V_{\emptyset ,s}$ on the invariant subspace $%
\mathcal{K}=\mathcal{H}_{\emptyset }$.
\end{enumerate}
\end{theorem}

\begin{proof}
Bearing in mind that $I_{k_{n}}>\dots >I_{k_{1}}>I_{k}$, and also $[\pi
_{k}(B),I_{k}]=0$, we represent the operators $F_{k}^{\Lambda }(\mathbf{b})$
determined by (2.5) in the form 
\begin{equation}
F_{k}^{\Lambda }(\mathbf{b})=I_{k_{n}}\pi _{k_{n}}(B_{n})\cdots I_{k_{1}}\pi
_{k_{1}}(B_{1})I_{k},\qquad k_{n}\leq \dots \leq k_{1}\leq k.
\end{equation}%
Hence, in accordance with (2.3), we obtain $\iota _{l}(a)F_{k}^{\Lambda }(%
\mathbf{b})=F_{k}^{\Lambda }(\mathbf{b})\iota _{l}(a)$, where $l=\max
\Lambda $, $a\in \mathfrak{A}_{l}$, whence we have the commutativity of $%
\kappa _{k}^{\Lambda }(\mathcal{B}_{k}^{\Lambda },\mathcal{B}_{k}^{\Lambda
})\subseteq \iota _{k}(\mathfrak{A}_{l})^{\prime }$.

The properties 1--4 are direct consequences of the corresponding conditions
of Definition 2 and can be directly verified for the compositions (2.5).
covariance holds for (2.5), as can be seen from condition 5, with respect to
the restrictions $V_{k,s}=V_{s}I_{k}$, which map the subspaces $\mathcal{H}%
_{k}=I_{k}\mathcal{H}$ (as follows from (2.4)) to $\mathcal{H}%
_{sk}:V_{s}I_{k}=I_{sk}V_{k}I_{k}$. At the same time, the operators $%
U_{s}=V_{\emptyset ,s}=V_{s}\mathrm{P}$, which leave $\mathcal{K}=\mathrm{P}%
\mathcal{H}$ invariant, are unitary, $U_{s}U_{s}^{\ast }=\mathrm{P}%
=U_{s}^{\ast }U_{s}$, by virtue of the fact, which follows from (2.4), that
the orthogonal projectors $U_{s}U_{s}^{\ast }$ commute with any of the
operators $a\in \mathfrak{A}_{\emptyset }=\mathcal{B}(\mathcal{K})$.

We now prove the self-consistency (0) of the family $\{\kappa _{k}^{\Lambda
}\}$. The condition 0 obviously holds, $\Lambda =\emptyset $, since $\kappa
_{j}^{\emptyset }=I_{j}$, and the empty sequence, of which $\mathcal{B}%
_{j}^{\emptyset }$ consists, induces the sequence $\mathbf{e}\in \mathcal{B}%
_{k}^{\mathrm{M}}$, for which $\kappa _{k}^{\mathrm{M}}(\mathbf{e},\mathbf{e}%
)=I_{k}$ by virtue of the normalization condition $P_{k_{i}}\leq I_{k_{i}}$
for any $\mathrm{M}\in \mathcal{F}$, $\mathrm{M}>k$.

Now suppose $\Lambda _{m}=\Lambda _{m-1}\cup j_{m}$, where $m\geq 1$, $%
j_{m}=\max \Lambda _{m}$, $\mathbf{b}_{m}=(\mathbf{b}_{m-1},B_{m})$ is the
sequence $\mathbf{b}_{m}\in \mathcal{B}_{j}^{\Lambda _{m}}$ determined by
the subsequence $\mathbf{b}_{m-1}\in \mathcal{B}_{j}^{\Lambda _{m-1}}$, $%
\Lambda _{m-1}=\Lambda _{m}\backslash j_{m}$, and $\mathrm{M}=\mathrm{M}%
_{n}\supseteq \Lambda _{m}$ is the finite subset $\mathrm{M}_{n}=\cup
_{i=1}^{n}k_{i}$, $n\geq m$, represented by the chain $k_{n}>\dots
>k_{1}>\emptyset $ of elements $k_{i}\in K$. We assume that the property 0
holds for $\Lambda =\Lambda _{m-1}$, $m\geq 1$, i.e., 
\begin{equation}
F_{j}^{\Lambda _{m-1}}(\mathbf{b}_{m-1})I_{k}=\pi _{k_{n}}(B_{n}^{\Lambda
_{m-1}})\pi _{k_{n-1}}(B_{n-1}^{\Lambda _{m-1}})\cdots \pi
_{k_{1}}(B_{1}^{\Lambda _{m-1}})I_{k},
\end{equation}%
where $(B_{1}^{\Lambda _{m-1}},\dots ,B_{n-1}^{\Lambda
_{m-1}},B_{n}^{\Lambda _{m-1}})\equiv \mathbf{b}_{m-1}^{\Lambda _{m-1}}$ is
determined by the induction $\mathbf{b}_{m-1}\mapsto \mathbf{b}%
_{m-1}^{\Lambda _{m-1}}$ from $\mathcal{B}^{\Lambda _{m-1}}$ to $\mathcal{B}%
^{\mathrm{M}_{n}}$: $B_{i}^{\Lambda _{m-1}}=B_{i}^{m-1}\times \mathrm{E}_{%
\bar{\Lambda}_{m-1}k_{i}}$, $i\leq n$. Note that $B_{n}=\mathrm{E}_{k_{n}}$,
since $\Lambda _{m-1}k_{n}=\emptyset $ ($\Lambda _{m-1}<j_{m}\leq k_{n}$),
and, taking into account conditions 0, 1, and 2, we obtain from (2.8), using
(2.7), 
\begin{equation}
F_{j}^{\Lambda _{m-1}}(\mathbf{b}_{m-1})I_{k}=\pi _{\Lambda
_{m-1}k_{n-1}}(B_{n-1}^{m-1})\cdots \pi _{\Lambda
_{m-1}k_{1}}(B_{1}^{m-1})I_{k}.
\end{equation}%
We now write $F_{j}^{\Lambda _{m}}(\mathbf{b}_{m})I_{k}=\pi
_{j_{m}}(B_{m})F_{j}^{\Lambda _{m-1}}(\mathbf{b}_{m-1})I_{k}$ in the form 
\begin{equation}
F_{j}^{\Lambda _{m}}(\mathbf{b}_{m})I_{k}=\pi _{j_{m}k_{n}}(B_{m_{n}})\pi
_{j_{m}k_{n-1}}(B_{m_{n-1}})\cdots \pi
_{j_{m}k_{1}}(B_{m_{1}})F_{j}^{\Lambda _{m-1}}(\mathbf{b}_{m-1})I_{k}
\end{equation}%
again using 0--2 and (2.7), where the elements $B_{mi}\in \mathcal{B}%
_{j_{m}k_{i}}$ are determined by the representation $B_{m}=\overset{n}{%
\underset{i=1}{\times }}B_{mi}$ which corresponds to the decomposition $%
j_{m}=\overset{n}{\underset{i=1}{\cup }}j_{m}k_{i}$. Substituting (2.9) in
(2.10) and bearing in mind that $j_{m}k_{i}\Join \Lambda _{m-1}k_{i^{\prime
}}$ for $i\leq i^{\prime }$, and also that $j_{m}k_{n}=\Lambda _{m}k_{n}$, $%
j_{m}k_{n-1}\cup \Lambda _{m-1}k_{n-1}=\Lambda _{m}k_{n-1},\dots
,j_{m}k_{1}\cup \Lambda _{m-1}k_{1}=\Lambda _{m}k_{1}$, we obtain by virtue
of the same conditions in the representation (2.7) 
\begin{equation}
F_{j}^{\Lambda _{m}}(\mathbf{b}_{m})I_{k}=\pi _{\Lambda
_{m}k_{n}}(B_{n}^{m})\pi _{\Lambda _{m}k_{n-1}}(B_{n-1}^{m})\cdots \pi
_{\Lambda _{m}k_{1}}(B_{1}^{m})I_{k},
\end{equation}%
where $B_{n}^{m}=B_{mn}$, $B_{n-1}^{m}=B_{mn-1}\times B_{n-1}^{m-1},\dots
,B_{1}^{m}=B_{m1}\times B_{1}^{m-1}$. The representation (2.7), like (2.9),
is equivalent to (2.8), and therefore 
\begin{equation}
F_{j}^{\Lambda _{m}}(\mathbf{b}_{m})I_{k}=\pi _{k_{n}}(B_{n}^{\Lambda
_{m}})\pi _{k_{n-1}}(B_{n-1}^{\Lambda _{m}})\cdots \pi
_{k_{1}}(B_{1}^{\Lambda _{m}})I_{k}.
\end{equation}%
Since (2.8) holds for $m=1$ ($\Lambda _{0}=\emptyset $), we obtain by
induction on $m$ the validity of (2.12), for any $m$, from which it follows
that property 0 holds for all $j<\Lambda \in \mathcal{F}$, which is what we
wanted to prove.
\end{proof}

\section{Existence and Reconstruction of QSP}

For every $l\in L$ we denote by $\mathfrak{A}_{l}=\underset{k\leq l}{\cap }%
\mathfrak{A}_{k}$, $T_{l}=\{t\in T:t\leq l\}$, and $\mathfrak{B}_{l}=%
\underset{\Lambda \leq l}{\cup }\mathcal{B}^{\Lambda }$ is the inductive
limit of the family $\mathcal{B}^{\Lambda }=\mathcal{B}_{\emptyset
}^{\Lambda }$, $\Lambda \in \mathcal{F}$, identified, as the set of
cylindric subsets $\times _{t\leq l}b\left( t\right) $, with the set of
functions $b:t\in T_{l}\mapsto b(t)\in \mathcal{B}_{t}$ equal to the unit $%
\mathrm{E}_{t}$ outside some finite $[\Lambda ]\subseteq T_{l}$. The sets $%
\mathfrak{B}_{l}$ are embedded into $\mathfrak{B}=\underset{\Lambda \in 
\mathcal{F}}{\cup }\mathcal{B}^{\Lambda }$ by the induction $b\in \mathfrak{B%
}_{l}\mapsto b_{l}\in \mathfrak{B}$, $b_{l}(t)=\mathrm{E}_{t}$ for $t\nleq l$%
, which extends $b=b_{l}|T_{l}$ to the function $b_{l}$ on the complete
factor-set $T$, and $\mathfrak{B}^{s}$, $s\in S$, will denote the set of
functions $b\in \mathfrak{B}$ that satisfy the condition $b(t)=\mathrm{E}%
_{t} $ if $t\notin sT$.

Let $\kappa :\mathfrak{B}\times \mathfrak{B}\rightarrow \mathcal{B}\left( 
\mathcal{K}\right) $ be the mapping with operator values uniquely determined
by the projective system $\{\kappa ^{\Lambda }\}$ of multikernels $\kappa
^{\Lambda }:\mathcal{B}^{\Lambda }\times \mathcal{B}^{\Lambda }\rightarrow 
\mathcal{B}\left( \mathcal{K}\right) $, 
\begin{equation}
\kappa ^{\Lambda }(\mathbf{b},\mathbf{b}^{\prime })=\kappa (\hat{\mathbf{b}},%
\hat{\mathbf{b}}^{\prime }),\qquad \mathbf{b},\mathbf{b}^{\prime }\in 
\mathcal{B}^{\Lambda },\quad \Lambda \in \mathcal{F}.
\end{equation}

The functional kernel $\kappa$ corresponding to the system of
operator-valued multikernels $\kappa^\Lambda(\mathbf{b},\mathbf{b}%
^{\prime})=\kappa_\emptyset^\Lambda (\mathbf{b},\mathbf{b}^{\prime})$
determined by (2.5) can be described by the following simple system of
axioms, which are equivalent to conditions 1--5 of Theorem 2 for the case $%
k=\emptyset$.

\begin{definition}
A $\mathcal{K}$-process with respect to $\mathfrak{A}$ in the wide sense
over $\mathcal{B}$ (or simply a QSP in the wide sense over $\mathcal{B}$ if $%
\mathcal{K}\simeq \mathbb{C}$) is a projective system of multikernels $%
\kappa =\{\kappa ^{\Lambda }:\Lambda \in \mathcal{F}\}$ described by the
mapping $\kappa :\mathfrak{B}\times \mathfrak{B}\rightarrow \mathcal{B}%
\left( \mathcal{K}\right) $ with values $\kappa (b_{l},b_{l}^{\prime })\in 
\mathfrak{A}_{l}^{\prime }$, $l\in L$, that satisfy the following axioms:

\begin{enumerate}
\item[$1^{0}$] $\sum_{i,i^{\prime }=1}^{m}(\zeta _{i}\mid \kappa
(b^{i},b^{i^{\prime }})\zeta _{i^{\prime }})\geq 0$, $\forall b^{i}\in 
\mathfrak{B}$, $\zeta _{i}\in \mathcal{K}$, $i\leq m=1,2,\dots $ .

\item[$2^{0}$] $\kappa (e,e)=1$, $e(t)=\mathrm{E}_{t}$, $\forall t\in T$.

\item[$3^0$] $\kappa(b_l,b_l)=\sum_{m=1}^m\kappa(b_l^m,b_l^m)$, if $%
b(t)=\sum_{m=1}^mb^m(t)$, $t\subseteq l$, $b^m(t)=b(t)$, $t\not\subseteq l$.

\item[$4^{0}$] $\kappa (b_{l}B,b_{l}^{\prime })=\kappa (b_{l},b_{l}^{\prime
}B)$ for every $B\in \mathcal{B}_{t}$, $t\subseteq l$, where $b,b^{\prime
}\in \mathfrak{B}_{l}$, $(bB)(t)=b(t)B$, and $(bB)(t)=b(t)$ for $%
t\not\subseteq l$. The process $\kappa $ in the wide sense is said to be $S$%
-covariant if

\item[$5^{0}$] $U_{s}^{\ast }\kappa (b,b^{\prime })U_{s}=\kappa
(b^{s},b^{\prime s})$, $\forall b,b^{\prime }\in \mathfrak{B}^{s}$, where $%
b^{s}(t)=b(st)^{s}$, $b^{\prime s}(t)=b^{\prime }(st)^{s}$, $t\in T$, $s\in
S $, with respect to the isometric representation $U=(U_{s})_{s\in S}$ of
the semigroup $S$ on $\mathcal{K}$. Every $\mathcal{H}$-process $(\mathcal{B}%
,\pi )$ with respect to the system $(\mathfrak{A},\iota )$ on $\mathcal{H}$
that determines a decomposition $\kappa ^{\Lambda }=\kappa _{\emptyset
}^{\Lambda }$ of the projective system of multikernels (3.1) in the form
(2.5) is called a realization of the process $\kappa $.
\end{enumerate}
\end{definition}

The existence of realizations for an arbitrary process $\kappa $, ensured
for the nonfunctional case of a single-point set $T$ by Naimark's
reconstruction theorem \cite{16}, is established by the following
(fundamental) theorem, which serves, in the case $\mathcal{K}\simeq \mathbb{C%
}$, as the noncommutative analog of Kolmogorov's reconstruction theorem
taking into account the causality relation on $T$.

\begin{theorem}
Let $\{\kappa ^{\Lambda }\}$ be the projective system (3.1) of mappings $%
\kappa ^{\Lambda }:\mathcal{B}^{\Lambda }\times \mathcal{B}^{\Lambda
}\rightarrow \mathcal{B}\left( \mathcal{K}\right) $ with values in $%
\mathfrak{A}_{l}^{\prime }$ for any $l\supseteq \max \Lambda $, determining
in accordance with the conditions $1^{0}-5^{0}$ an $S$-covariant $\mathfrak{A%
}$ process with respect to $\mathcal{K}$ over $\mathcal{B}$ in the wide
sense. Then there exists a canonical system $\pi ^{\ast }$ of $\sigma $%
-representations of the family $\mathcal{B}$ which is described by
conditions 1--5 of Definition 2 with respect to a representation system $(%
\mathfrak{A},\iota )$ on $\mathcal{H}$ and an isometric $S$-representation $%
V $ as $S$-covariant $\mathcal{H}$-process that determines the decomposition 
$\kappa ^{\Lambda }=\kappa _{\emptyset }^{\Lambda }$ (2.5) and satisfies the
condition $\pi _{k}^{\ast }(\mathrm{E}_{k})=\underset{l\supseteq k^{\prime }}%
{\vee }E_{l}\equiv P_{k}^{\ast }$ for all $k\in K$, where $E_{l}=\underset{%
k\subseteq l}{\wedge }P_{k}^{\ast }$, $l\in L$. If in addition $\mathfrak{A}%
_{k}=\underset{l\subseteq k}{\vee }\mathfrak{\mathfrak{A}_{l}}$, $k\neq
\emptyset $, and for any $b\in \mathfrak{B}$, $\zeta \in \mathcal{K}$ 
\begin{equation}
\underset{l\downarrow \emptyset }{\lim }\underset{\{\zeta _{i},b^{i}\}}{\sup 
}\left\vert \sum_{i}(\zeta _{i}|\kappa (b^{i},b)-\kappa (b^{i},e)\kappa
(e,b)|\zeta )\right\vert =0,
\end{equation}%
where $\zeta _{i}\in \mathcal{K}$, $b^{i}\in \mathfrak{B}_{l}$, $i\in 
\mathbf{N}$, $\sum (\zeta _{i}|\kappa (b^{i},b^{i^{\prime }})\zeta
_{i^{\prime }})\leq 1$, then there exists a realizing $\mathcal{H}$-process $%
(\mathcal{B},\pi )$ with respect to $\ast $-representations $\iota ^{\ast }$
of the family $\mathfrak{A}$ that satisfy the conditions $\iota _{k}^{\ast
}(1)=\underset{l\supseteq k}{\wedge }E_{l}$ for all $k\in K$. Moreover, if
for each sequence $B_{i}=\underset{t\subseteq l_{i}}{\times }b(t)b^{\prime
}\left( t\right) $, $i=1,\dots ,n$, corresponding to a chain $l_{1}\leq
\dots \leq l_{n}$ the kernel $\kappa $ is determined by a regression 
\begin{equation}
\kappa (b,b^{\prime })=\rho _{l_{1}}(\pi _{l_{1}}(B_{1})\theta
_{l_{1}l_{2}}(\pi _{l_{2}}(B_{2})\cdots \theta _{l_{n-1}l_{n}}(\pi
_{l_{n}}(B_{n}))\cdots )),
\end{equation}%
where $\rho _{l}:\mathcal{A}_{l}\rightarrow \mathcal{B}\left( \mathcal{K}%
\right) $, $l\in L$ and $\theta _{ll^{\prime }}:\mathcal{A}_{l^{\prime
}}\rightarrow \mathcal{A}_{l}$ are normal mappings determined for all $l\leq
l^{\prime }\in L$ on the W*-algebra $\mathcal{A}=\pi _{l}(\mathcal{B}%
_{l})\vee \iota _{l}(\mathfrak{A}_{l})$ generated by representations $\pi
_{l}:\mathcal{B}_{l}\rightarrow \mathcal{P}(\mathcal{H}_{l})$, $\iota _{l}:%
\mathfrak{A}_{l}\rightarrow \mathfrak{\pi }_{l}(\mathcal{B}_{l})^{\prime }$
on some Hilbert spaces $\mathcal{H}_{l}$, $l\in L$, then the canonical
realization $\pi ^{\ast }$ has the conditional Markov property 
\begin{equation}
E_{l}\pi _{l^{\prime }}^{\ast }(\mathcal{B}_{l^{\prime }})E_{l}\subseteq \pi
_{l}^{\ast }(\mathcal{B}_{l})\vee \iota _{l}^{\ast }(\mathfrak{A}%
_{l}),\qquad \forall l\leq l^{\prime }\in L.
\end{equation}
\end{theorem}

We divide the proof into five stages.

1. \textbf{Reconstruction of the space} $\mathcal{H}$. We consider the $%
\mathcal{K}$-hull of the set $\mathfrak{B}$, which is defined as the linear
space of formal sums $\zeta _{i}b^{i}$ ($=\sum_{i}\zeta _{i}b^{i}$) with
finite number of nonvanishing coefficients $\zeta _{i}\in \mathcal{K}$ with
respect to the componentwise defined operations 
\begin{equation*}
\zeta _{i}b^{i}+\zeta _{i}^{\prime }b^{i}=(\zeta _{i}+\zeta _{i}^{\prime
})b^{i},\quad c(\zeta _{i}b^{i})=(c\zeta _{i})b^{i},\quad c\in \mathbb{C}.
\end{equation*}%
Let $\mathcal{E}$ be the factorisation of this $\mathcal{K}$-hull with
respect to the non-negative definite (by virtue of $1^{0}$) Hermitian form 
\begin{equation}
(\zeta _{i}b^{i}\mid \zeta _{i^{\prime }}b^{i^{\prime }})=(\zeta _{i}\mid
\kappa (b^{i},b^{i^{\prime }})\zeta _{i^{\prime }}=)\equiv \Vert \zeta
_{i}b^{i}\Vert ^{2},
\end{equation}%
and $\mathcal{H}=\bar{\mathcal{E}}$ be the completion of $\mathcal{E}$ with
respect to the norm $\Vert \cdot \Vert $. We denote by $|\zeta b)=\{\zeta
_{i}b^{i}:\zeta _{i}b^{i}-\zeta b\in \mathcal{N}\}$ the equivalence classes
of elements $\zeta b$ generating $\mathcal{E}$, determined by the kernel $%
\mathcal{N}=\{\zeta _{i}b^{i}:\Vert \zeta _{i}b^{i}\Vert =0\}$ of the form
(3.1), retaining it as notation of the scalar product of the bra $(\zeta
_{i}b^{i}|$ and ket $|\zeta _{i}b^{i})=(\zeta _{i}b^{i}|^{\ast }$ vectors.
The space $\mathcal{K}$ can be embedded into $\mathcal{H}$ by identification
of each $\zeta \in \mathcal{K}$ and vector $|\zeta e)$ isometric to it by
virtue of $2^{0}$. Let $\mathcal{E}_{j}\subseteq \mathcal{E}$ for each $j\in
J$ be the subspace generated by the vectors $|\zeta b)$, $\zeta \in \mathcal{%
K}$, $b\in \mathfrak{B}_{j}=\underset{\Lambda \leq j}{\cup }\mathcal{B}%
^{\Lambda }$, and $E_{j}$ be the orthogonal projector onto the completions $%
\mathcal{H}_{j}=\bar{\mathcal{E}}_{j}$, which obviously satisfies the
condition $E_{j}=\underset{k\subseteq j}{\vee }E_{k}$, where $k\in K$, $j\in
J$, $E_{\emptyset }=\mathrm{P}$ is the orthogonal projector onto the
subspace $\mathcal{K}\subseteq \mathcal{H}$. Let $P_{j}^{\ast }=\underset{%
l\supseteq j}{\vee }E_{l}$, $I_{j}^{\ast }=\underset{l\supseteq j}{\wedge }%
E_{l}$ be the orthogonal projectors onto the completions $\mathcal{K}%
_{j}^{\ast }=\bar{\mathcal{E}}_{j[}$, $\mathcal{H}_{j}^{\ast }=\bar{\mathcal{%
E}}_{j]}$ of the subspaces $\mathcal{E}_{j[}=\underset{l\supseteq j}{\cup }%
\mathcal{E}_{l}$, $\mathcal{E}_{j]}=\underset{l\supseteq j}{\cap }\mathcal{E}%
_{l^{\prime }}$ generated by the vectors $|\zeta b)$, $\zeta \in \mathcal{K}$%
, respectively, for $b\in \mathfrak{B}_{j[}:=\underset{l\supseteq j}{\cup }%
\mathfrak{B}_{l}$ and $b\in \mathfrak{B}_{j]}:=\underset{l\supseteq j}{\cap }%
\mathfrak{B}_{l}$. Because $\underset{k\subseteq l}{\cap }\mathfrak{B}_{k[}=%
\mathfrak{B}_{l}$ for any $l\in L$ and $\mathfrak{B}_{\emptyset \lbrack }=%
\mathfrak{B}$, we obtain condition 0 for $P_{k}=P_{k}^{\ast }$ and $%
I_{k}=E_{k}$, $k\in K$: $\underset{k\supseteq j}{\wedge }P_{k}^{\ast }=E_{l}$%
, $P_{\emptyset }^{\ast }=\mathrm{I}$. Since $\mathfrak{B}_{k[}\subseteq 
\mathfrak{B}_{l}$ for $k\subseteq l$, we also have $\underset{k\subseteq l}{%
\cup }\mathfrak{B}_{k]}=\mathfrak{B}_{l}$ and $\underset{k\subseteq l}{\vee }%
I_{k}^{\ast }=E_{l}$ for $l\in L$. We show that $I_{\emptyset }^{\ast }=%
\mathrm{P}$, i.e., $\underset{l\in L}{\wedge }E_{l}=\underset{l\downarrow
\emptyset }{\lim }E_{l}=\mathrm{P}$ (by virtue of the monotonicity of $E_{t}$%
) when the regularity condition (3.2) is satisfied. Since $E_{l}\geq \mathrm{%
P}$, 
\begin{equation*}
(\zeta b|\left( E_{l}-\mathrm{P}\right) |\zeta b)=(\zeta _{\cdot }b^{\cdot
}|E_{l}|\zeta _{\cdot }b^{\cdot })=\Vert E_{l}\zeta _{\cdot }b^{\cdot }\Vert
^{2},
\end{equation*}%
where $|\zeta _{\cdot }b^{\cdot })=|\zeta b)-\mathrm{P}|\zeta b)=|\zeta
b)-|\kappa (e,b)\zeta e)$. By virtue of the fact that the family $|\zeta
_{i}b^{i})$, $\zeta _{i}\in \mathcal{K}$, $b^{i}\in \mathfrak{B}_{l}$, is
dense in the subspace $\mathcal{H}_{l}=E_{l}\mathcal{H}$ and by the
definition of the norm (3.5), we have%
\begin{eqnarray*}
\Vert E_{l}\zeta _{\cdot }b^{\cdot }\Vert ^{2} &=&\sup \{|(\zeta
_{i}b^{i}|\zeta _{\cdot }b^{\cdot })|^{2}:\;\Vert \zeta _{i}b^{i}\Vert
^{2}\leq 1\} \\
&=&\sup \{|(\zeta _{i}|\kappa (b^{i},b^{\cdot })\zeta _{\cdot })|^{2}:(\zeta
_{i}|\kappa (b^{i},b^{i^{\prime }})\zeta _{i^{\prime }})\leq 1\} \\
&=&\underset{\{b^{i},\zeta _{i}\}}{\sup }|(\zeta _{i}|\kappa
(b^{i},b)-\kappa (b^{i},e)\kappa (e,b)|\zeta )|^{2}\downarrow 0
\end{eqnarray*}%
for $l\downarrow \emptyset $, $\zeta _{i}\in \mathcal{K}$, $b^{i}\in 
\mathfrak{B}_{l}$, $(\zeta _{i}|\kappa (b^{i},b^{i^{\prime }})\zeta
_{i^{\prime }})\leq 1$. Thus, $E_{l}\downarrow \mathrm{P}$ for $l\downarrow
\emptyset $, $l\in L$ in the strong operator topology.

2. \textbf{Reconstruction of the representation} $\iota $. For any $a\in 
\mathfrak{A}_{j}=\underset{k\subseteq j}{\cap }\mathfrak{A}_{k}$, we set 
\begin{equation}
\lambda _{j}(a)|\zeta _{i}\hat{b}^{i})=|a\zeta _{i}\hat{b}^{i}),\qquad \zeta
_{i}\in \mathcal{K},\quad b^{i}\in \mathcal{B}^{\Lambda },\quad \Lambda \leq
j\in J.
\end{equation}%
The operator on $\mathcal{E}_{j}$ defined in this manner is linear and
identical to $\lambda _{j^{\prime }}(a)|\mathcal{E}_{j}$ for $j\leq
j^{\prime }$, and this definition is correct by virtue of the inequality 
\begin{equation*}
\Vert a\zeta _{i}\hat{b}^{i}\Vert ^{2}=(a\zeta _{i}|\kappa ^{\Lambda
}(b^{i},b^{i^{\prime }})a\zeta _{i^{\prime }})\leq \Vert a\Vert ^{2}\Vert
\zeta _{i}\hat{b}^{i}\Vert ^{2},
\end{equation*}%
which follows from the condition of commutativity $\kappa ^{\Lambda
}(b,b^{\prime })\in \mathfrak{A}_{j}^{\prime }$ with $a\in \mathfrak{A}_{j}$
for $b^{\prime }\in \mathcal{B}^{\Lambda }$, $\Lambda \leq j$. The mapping $%
a\mapsto \lambda _{j}(a)$ is obviously linear and multiplicative with
respect to $a\in \mathfrak{A}_{j}$, and $\lambda _{j}(a^{\ast })=\lambda
_{j}(a)^{\ast }|\mathcal{E}_{j}$, and $\lambda _{j}(1)=\mathrm{I}|\mathcal{E}%
_{j}$. Denoting by $\iota _{j}(a^{\ast })=\lambda _{j}(a)^{\ast }$ the
continuous extension of the bounded operators $\lambda _{j}(a)^{\ast }$ to
the whole of $\mathcal{H}$ determined by the condition $\iota _{j}(a^{\ast
})|\mathcal{E}_{j}^{\bot }=0$ on the orthogonal complement $\mathcal{E}%
_{j}^{\bot }=\mathcal{H}_{j}^{\bot }$, we obtain the family of $\ast $%
-representations $\iota _{j}:\mathfrak{A}_{j}\rightarrow \mathcal{B}\left( 
\mathcal{H}\right) $, which satisfy the condition (2.3) with respect to $%
\iota _{j}(1)=I_{j}=E_{j}$, $j\in J$, and $\iota _{\emptyset }$ is obviously
the identity representation of the C*-algebra $\mathfrak{A}_{\emptyset }=%
\mathcal{B}\left( \mathcal{K}\right) $ on the subspace $\mathcal{K}\subseteq 
\mathcal{H}$.

Denoting $\lambda _{j}^{l}(a)=\lambda _{l}(a)|\mathcal{E}_{j]}$ for fixed $%
j\leq l\in L$ and taking into account the invariance with respect to $%
\lambda _{l}(a)$ of the subspaces $\mathcal{E}_{l^{\prime }}$, $l^{\prime
}\leq l$, whose intersection for $l^{\prime }\supseteq j$ is $\mathcal{E}%
_{j]}$, we obtain a family of self-consistent $\ast $-representations $%
\lambda _{j}^{l}:\mathfrak{A}_{l}\rightarrow \mathcal{B}(\mathcal{E}_{j]})$, 
$\lambda _{j}^{l^{\prime }}|\mathfrak{A}_{l}=\lambda _{j}^{l}$ for any $%
j\leq l^{\prime }\leq l$, that determine the projective limit $\lambda
_{j}^{\ast }:\mathfrak{A}_{j}^{\ast }\rightarrow \mathcal{B}(\mathcal{E}%
_{j]})$ by $\lambda _{j}^{\ast }|\mathfrak{A}_{l}=\lambda _{j}^{l}$ for $%
l\supseteq j$, this representing the $\ast $-algebra $\mathfrak{A}_{j}^{\ast
}=\underset{l\supseteq j}{\vee }\mathfrak{A}_{l}$ on $\mathcal{E}_{j]}$ with 
$\lambda _{j}^{\ast }(1)=\mathrm{I}|\mathcal{E}_{j]}$. The continuous
extension $\iota _{j}^{\ast }(a^{\ast })=\lambda _{j}^{\ast }(a)^{\ast }$,
equal to zero on $\mathcal{E}_{j]}^{\bot }$, also satisfies the condition
(2.3) with respect to $\iota _{j}^{\ast }(1)=I_{j}^{\ast }=\underset{%
l\supseteq j}{\wedge }E_{l}$, determining self-consistent $\ast $%
-representations $\iota _{j}^{\ast }:\mathfrak{A}_{j}^{\ast }\rightarrow 
\mathcal{B}\left( \mathfrak{\mathcal{H}}\right) $ that are identical on the
subspace $\mathcal{H}_{j}\subseteq \mathcal{H}_{j}^{\ast }$ to the
representations $\iota _{j}$ restricted to the $\ast $-subalgebras $%
\mathfrak{A}_{j}^{\ast }\subseteq \mathfrak{A}_{j}$ (for $j\in L\;\;\;%
\mathfrak{A}_{j}^{\ast }=\mathfrak{A}_{j}$, $\mathcal{H}_{j}^{\ast }=%
\mathcal{H}_{j}$, and $\iota _{j}^{\ast }=\iota _{j}$). If $\mathfrak{A}%
_{k}^{\ast }=\mathfrak{A}_{k}$ for any $k\in K$, $k\neq \emptyset $, then as
the required family of representations $\mathfrak{A}=(\mathfrak{\mathfrak{A}}%
_{k})_{k\in K}$ we can also choose $\iota ^{\ast }=(\iota _{k}^{\ast
})_{k\in K}$.

3. \textbf{Reconstruction of the representation} $\pi $. Identifying the
elements $b\in \mathfrak{B}_{l}$ with functions $b\in \mathfrak{B}$ equal to
unity for $t\not\leq l$, we set for each $B\in \mathcal{B}_{j}\subseteq 
\underset{t\in \lbrack j]}{\otimes }\mathcal{B}_{t}$, $j\subseteq l\in L$, 
\begin{equation}
\rho _{l}^{j}(B)|\zeta _{i}b^{i})=|\zeta _{i}b^{i}B),\qquad \forall \zeta
_{i}\in \mathcal{K},\quad b^{i}\in \mathfrak{B}_{l},
\end{equation}%
where the product $b^{i}B\in \mathfrak{B}$ is the function equal to $%
(b^{i}B)(t)=b^{i}(t)b(t)$ when $t\in \lbrack j]$ and $(b^{i}B)(t)=b^{i}(t)$
when $t\notin \lbrack j]$ for $B=\underset{t\in \lbrack j]}{\times }b(t)$.
Since $(b^{i}B)B=b^{i}B$, the operator $\rho _{l}^{j}(B)$ is idempotent on $%
\mathcal{E}_{l}$ and satisfies the condition 
\begin{equation*}
(\zeta _{i}b^{i}B|\zeta _{i^{\prime }}b^{i^{\prime }})=(\zeta
_{i}b^{i}|\zeta _{i^{\prime }}b^{i^{\prime }}B)=\Vert \zeta _{i}b^{i}B\Vert
^{2},
\end{equation*}%
which defines it correctly on $\mathcal{E}_{l}$ as a linear symmetric
projector: 
\begin{equation*}
\rho _{l}^{j}(B)^{\ast }|\mathcal{E}_{l}=\rho _{l}^{j}(B)=\rho
_{l}^{j}(B)^{2}|\mathcal{E}_{l},
\end{equation*}%
identical to $\rho _{l^{\prime }}^{j}(B)|\mathcal{E}_{l}$ for any $l^{\prime
}\in L$ such that $l\leq l^{\prime }$, $j\subseteq l^{\prime }$. The mapping 
$B\mapsto \rho _{l}^{j}(B)$ satisfies on $\mathcal{E}_{l}$ the following
conditions for $j\subseteq l\in L$, these following directly from the
definition (3.7) for all $\zeta _{i}\in \mathcal{K}$, $b^{i}\in \mathfrak{B}%
_{l}$:

\begin{enumerate}
\item[0)] $\rho _{l}^{j}(\mathrm{E}_{j})|\zeta _{i}b^{i})=|\zeta _{i}b^{i}%
\mathrm{E}_{j})=|\zeta _{i}b^{i})\qquad (\mathrm{E}_{j}=\underset{t\in
\lbrack j]}{\times }\mathrm{E}_{t})$;

\item[1)] $j\sim j^{\prime}\Rightarrow\rho_l^{j\cup
j^{\prime}}(BB^{\prime})|\zeta_ib^i)=
\rho_l^j(B)\rho_l^{j^{\prime}}(B^{\prime})|\zeta_ib^i)$;

\item[2)] $j\Join j^{\prime}\Rightarrow\rho_l^{j\cup j^{\prime}}(B\times
B^{\prime})|\zeta_ib^i)=
\rho_l^j(B)\rho_l^{j^{\prime}}(B^{\prime})|\zeta_ib^i)$, where $l\supseteq
j\cup j^{\prime}$, $B\in\mathcal{B}_j$, $B^{\prime}\in\mathcal{B}%
_{j^{\prime}}$;

\item[3)] for any decomposition $\mathrm{E}_{j}=\sum_{m=1}^{\infty }B^{m}$
it follows from $3^{0}$ that 
\begin{equation*}
(\zeta _{i}b^{i}|\zeta _{i^{\prime }}b^{i^{\prime }})=\sum_{m=1}^{\infty
}(\zeta _{i}b^{i}B^{m}|\zeta _{i^{\prime }}b^{i^{\prime
}}B^{m})=\sum_{m=1}^{\infty }(\zeta _{i}b^{i}|\rho _{l}^{j}(B^{m})|\zeta
_{i^{\prime }}b^{i^{\prime }});
\end{equation*}

\item[4)] $[\iota _{j}(a),\rho _{l}^{j}(B)]|\zeta _{i}b^{i})=0$, $\forall
a\in \mathfrak{A}_{j}$, $B\in \mathcal{B}_{j}$ (and similarly for $\iota
_{j}^{\ast }(a)$, $a\in \mathfrak{A}_{j}^{\ast }$), since the projectors $%
\rho _{l}^{j}(B)$ commute with both $\iota _{j}(\mathfrak{A}_{j})$ and $%
\iota _{j}^{\ast }(\mathfrak{A}_{j}^{\ast })$ on the invariant subspaces $%
\mathcal{E}_{j}$ and $\mathcal{E}_{j]}\subseteq \mathcal{E}$. The symmetric
projectors $\rho ^{j}(B)$, which are uniquely determined on $\mathcal{E}%
_{j[}=\underset{l\supseteq j}{\cup }\mathcal{E}_{l}$ by the condition $\rho
^{j}(B)|\mathcal{E}_{l}=\rho _{l}^{j}(B)$, $l\supseteq j$, also have
properties 0--4 and, therefore, can be extended by continuity to the Hilbert
subspace $\mathcal{K}_{j}^{\ast }=\bar{\mathcal{E}}_{j[}$ to Hermitian
projectors $\pi _{j}^{\ast }(B)=\rho ^{j}(B)^{\ast }$, $\pi _{j}^{\ast }(B)|%
\mathcal{E}_{j[}^{\bot }=0$, satisfying conditions 0--4 of Definition 2.
Noting that $\pi _{j}^{\ast }(\mathrm{E}_{j})$ is an orthogonal projector $%
P_{j}^{\ast }$ onto $\mathcal{K}_{j}^{\ast }$, we obtain $P_{j}^{\ast
}=P_{j^{\prime }}^{\ast }$, if $j\sim j^{\prime }$, $P_{j}^{\ast }\wedge
P_{j^{\prime }}^{\ast }=P_{j\cup j^{\prime }}^{\ast }$, if $j\Join j^{\prime
}$, and $P_{\emptyset }^{\ast }=\mathrm{I}$.
\end{enumerate}

4. \textbf{Reconstruction of the representation} $V$. For any $s\in S$, we
set 
\begin{equation}
U^{s}|\zeta _{i}b^{i})=|U_{s}^{\ast }\zeta _{i}b^{is}),\qquad \forall \zeta
_{i}\in \mathcal{K},\quad b^{i}\in \mathfrak{B}^{s},
\end{equation}%
where $b^{s}(t)=b(st)^{s}$ is determined by functions $b\in \mathfrak{B}^{s}$
equal to identity $b(t)=\mathrm{E}_{t}$ for $t\notin \lbrack s\mathrm{T}]$.
The operator $U^{s}$ is correctly defined on the pre-Hilbert subspace $%
\mathcal{E}^{s}$ generated by the vectors $|\zeta b)$ for $\zeta \in 
\mathcal{K}$, $b\in \mathfrak{B}^{s}$, by virtue of condition $5^{0}$, which
reduces to the condition that the operator be isometric: 
\begin{align*}
(\zeta _{i}b^{i}|U^{s\ast }U^{s}|\zeta _{i^{\prime }}b^{i^{\prime }})&
=(U_{s}^{\ast }\zeta _{i}|\kappa (b^{is},b^{i^{\prime }s})U_{s}^{\ast }\zeta
_{i^{\prime }})=(\zeta _{i}|U_{s}\kappa (b^{is},b^{i^{\prime }s})U_{s}^{\ast
}|\zeta _{i^{\prime }})= \\
& (\zeta _{i}b^{i}|\zeta _{i^{\prime }}b^{i^{\prime }}),
\end{align*}%
where we have used the fact that $U_{s}U_{s}^{\ast }\zeta =\zeta $ for $%
\zeta \in \mathcal{K}$. By virtue of the assumed surjectivity of the
mappings $B\in \mathcal{B}_{sj}\mapsto B^{s}\in \mathcal{B}_{j}$, the
operators $U^{s}$ map every subspace $\mathcal{E}_{sj}^{s}=\mathcal{E}%
^{s}\cap \mathcal{E}_{sj}$ isometrically on $\mathcal{E}_{j}$, i.e., are
unitary on these subspaces and are made self-consistent by the condition $%
U^{s}U^{s^{\prime }}|\mathcal{E}^{ss^{\prime }}=U^{ss^{\prime }}$, $\forall
s,s^{\prime }\in S$. Using the definitions (3.6) and (3.7) of the
representations $\lambda _{j}$ and $\rho ^{j}$, we obtain for $a\in 
\mathfrak{A}_{sj}$, $b^{i}\in \mathfrak{B}_{sj}^{s}$ 
\begin{equation*}
\lambda _{j}(a^{s})U^{s}|\zeta _{i}b^{i})=|a^{s}\zeta
_{i}b^{is})=U^{s}\lambda _{sj}(a)|\zeta _{i}b^{i})
\end{equation*}%
and similarly for $\lambda _{j}^{\ast }$ when $a\in \mathfrak{A}_{sj}^{\ast
} $, $b^{i}\in \mathfrak{B}_{sj]}^{s}$, and also 
\begin{equation*}
\rho ^{j}(B^{s})U^{s}|\zeta _{i}b^{i})=|\zeta _{i}^{s}b^{is}B^{s})=U^{s}\rho
^{sj}(B)|\zeta _{i}b^{i})
\end{equation*}%
for $b^{i}\in \mathfrak{B}_{sj[}^{s}$, $B\in \mathcal{B}_{sj}$. Thus, the
operators $U^{s}$ intertwine the representations $\iota _{j}^{s}(a)=\iota
_{j}(a^{s})$ on $\mathcal{E}_{j}$ with $\iota _{sj}(a)$ on $\mathcal{E}%
_{sj}^{s}\subseteq \mathcal{E}_{sj}$, $\iota _{j}^{\ast s}(a)=\iota
_{j}^{\ast }(a^{s})$ on $\mathcal{E}_{j]}$ with $\iota _{sj}^{\ast }(a)$ on $%
\mathcal{E}_{sj]}^{s}\subseteq \mathcal{E}_{sj]}$, and $\pi _{j}^{\ast
s}(B)=\pi _{j}^{\ast }(B^{s})$ on $\mathcal{E}_{j[}$ with $\pi _{sj}^{\ast
}(B)$ on $\mathcal{E}_{sj[}$. Denoting by $V_{s}$ the inverse unitary
operators $\mathcal{E}\rightarrow \mathcal{E}^{s}$, extended by continuity
to linear isometries $V_{s}:\mathcal{H}\rightarrow \mathcal{H}$ that map $%
\mathcal{H}$ onto $\mathcal{H}^{s}=\bar{\mathcal{E}}^{s}$, we obtain from
this, bearing in mind that $V_{s}^{\ast }|\mathcal{E}^{s}=U^{s}$, the
required representation $V=(V_{s})_{s\in S}$, which determines for the
constructed representations $\iota =(\iota _{k})_{k\in K}$ or for $\iota
^{\ast }$ and $\pi ^{\ast }=(\pi _{k}^{\ast })_{k\in K}$ the covariance
conditions (2.5) and condition 5 of Definition 2 and is identical to $%
U=(U_{s})_{s\in S}$ on the subspace $\mathcal{H}_{\emptyset }=\mathcal{K}$.

5. \textbf{Decomposition of the multikernel} $\kappa $. Since any pair $%
b,b^{\prime }\in \mathfrak{B}$ of functions $t\mapsto b(t)$, $b(t^{\prime
})\in \mathcal{B}_{t}$ has common set $\Lambda \in \mathcal{F}$, for which $%
b(t)=E_{t}=b^{\prime }(t)$ for $t\notin \lbrack \Lambda ]$, there exists a
finite chain $\Lambda _{n}=\{k_{1},\dots ,k_{n}\}$, $k_{n}>\dots
>k_{1}>\emptyset $, $\overset{n}{\underset{i=1}{\cup }}k_{i}=\Lambda $, and
a pair of sequences $\mathbf{b}=(B_{1},\dots ,B_{n})$, $\mathbf{b}^{\prime
}=(B_{1}^{\prime },\dots ,B_{n}^{\prime })$ determined by functions on $%
[\Lambda ]$ having continuations $b=\hat{b}$ and $b^{\prime }=\hat{b}%
^{\prime }$ on $T$ that represent the corresponding value of the kernel $%
\kappa $ in the form 
\begin{equation}
\kappa (b,b^{\prime })=\kappa (\hat{b},\hat{b}^{\prime }).
\end{equation}%
We define the chronologically ordered product $F^{\Lambda _{n}}(\mathbf{b})$
by means of the constructed representation $\pi $, acting successively on an
arbitrary vector $\zeta \in \mathcal{K}$: 
\begin{equation*}
F^{\Lambda _{n}}(\mathbf{b})\zeta =\pi _{k_{n}}(B_{n})\cdots \pi
_{k_{1}}(B_{1})|\zeta e)=|\zeta b^{n}),
\end{equation*}%
where $b^{n}\in \mathfrak{B}$ is a function of $t\in \lbrack T]\mapsto 
\mathcal{B}_{t}$ determined recursively as follows: $b^{n}(t)=b^{n-1}(t)$
for $t\notin \lbrack k_{n}]$, $b^{n}(t)=b_{n}(t)$, $t\in \lbrack k_{n}]$
with the initial function $b^{0}\left( t\right) =e(t)=\mathrm{E}_{t}$ and
the functions $b_{n}$ defining the representation $B_{n}=\underset{t\in
\lbrack k_{n}]}{\times }b_{n}(t)$, $n=1,2,\dots $, this being obviously
identical to $b=\hat{b}$. Therefore, we obtain from (3.1) and (3.5) 
\begin{equation*}
(\zeta |F^{\Lambda _{n}}(b)^{\ast }F^{\Lambda _{n}}(b)|\zeta ^{\prime
})=(\zeta b|\zeta ^{\prime }b^{\prime })=(\zeta |\kappa ^{\Lambda
_{n}}(b,b^{\prime })|\zeta ^{\prime }).
\end{equation*}%
Thus, the obtained family $\pi ^{\ast }$, which represents $\mathcal{B}$ on $%
\mathcal{H}$, determines the decomposition (2.5) and satisfies all the
conditions of Definition 2 with respect to the family $\iota $ and the
representations of $\mathfrak{A}$, and also with respect to $\iota ^{\ast }$
for $\mathfrak{A}^{\ast }=\mathfrak{A}$. This proves the existence of a
realization of the process defined in the wide sense.

We now show that when the conditions of the regression (3.3) are satisfied,
the constructed representations $\pi _{l}^{\ast }$, $\iota _{l}^{\ast }$ ($%
=\iota _{l}$) for $l\in L$ have the property (3.4) with respect to the
orthogonal projectors $\iota _{l}^{\ast }(1)=E_{l}=\pi _{l}^{\ast }(\mathrm{E%
}_{l})$. Using the definition (3.7) of the representations $\pi _{l}^{\ast
}(B)|\mathcal{E}_{l[}=\rho _{l}^{l}(B)$, $B\in \mathcal{B}_{l}$, we have 
\begin{equation*}
E_{l}\pi _{l^{\prime }}^{\ast }(B)|\zeta b)=\lambda _{l}\circ \nu
_{ll^{\prime }}(B)|\zeta b),\quad l\leq l^{\prime },\quad \zeta \in \mathcal{%
K},\quad b\in \mathfrak{B}_{l},
\end{equation*}%
where $\nu _{ll^{\prime }}=\theta _{ll^{\prime }}\circ \pi _{l^{\prime }}$, $%
B\in \mathcal{B}_{l^{\prime }}$, and $\lambda _{l}:\mathcal{A}%
_{l}\rightarrow \mathcal{B}\left( \mathcal{H}\right) $ is the
W*-representation uniquely determined by the condition $\lambda
_{l}(AP)=\iota _{l}^{\ast }(a)\pi _{l}^{\ast }(B)$ for all $A=\iota _{l}(a)$%
, $a\in \mathfrak{A}_{l}$ and $P=\pi _{l}(B)$, $B\in \mathcal{B}_{l}$. Since
obviously $\lambda _{l}(AP)A^{\prime }=A^{\prime }\lambda _{l}(AP)$ for any $%
A\in \iota _{l}(\mathfrak{A}_{l})$ and $P\in \pi _{l}(\mathcal{B}_{l})$, $%
A^{\prime }\in \pi _{l}^{\ast }(\mathcal{B}_{l})^{\prime }\cap \iota
_{l}^{\ast }(\mathfrak{A}_{l})^{\prime }$, we obtain $A^{\prime }E_{l}\pi
_{l}^{\ast }(B)E_{l}=E_{l}\pi _{l}^{\ast }(B)E_{l}A^{\prime }$, i.e., the
condition (3.4). Thus, the reconstruction theorem, Theorem 3, is proved.

\section{Equivalence and Interpretation of QSP}

1. Let $\mathfrak{A}_{l}=\underset{k\subseteq l}{\cap }\mathfrak{A}_{k}$ be
a nonincreasing family of $\ast $-subalgebras $\mathfrak{A}_{l}\subseteq 
\mathcal{B}(\mathcal{K})$, $\mathfrak{A}_{l}\supseteq \mathfrak{A}%
_{l^{\prime }}$, if $l\leq l^{\prime }$, indexed by the maximal subsets $%
l\subset T$ of pairwise nonanticipatory elements $t\in T$, and $\mathcal{B}%
_{l}=\underset{t\in l}{\times }\mathcal{B}_{t}$ be Boolean semirings of
cylindrical subsets $B=\underset{t\in l}{\times }b(t)$, determined by
functions $b:t\mapsto b(t)\in \mathcal{B}_{t}$ equal to $\mathrm{E}_{t}$
outside a certain finite subset of equivalence classes $t=[t]\subseteq l$.
For every $\mathcal{H}$-process $(\mathcal{B}_{k},\pi _{k})_{k\in K}$,
defined with respect to the system $(\mathfrak{A}_{k},\iota _{k})_{k\in K}$,
we denote by $\iota _{l},\pi _{l}$ for each $l\in L$, the $\ast $%
-representation $\mathfrak{A}_{l}\rightarrow \mathcal{B}\left( \mathcal{H}%
\right) $ and $\sigma $-representation $\mathcal{B}_{l}\rightarrow \mathcal{P%
}\left( \mathcal{H}\right) $ uniquely determined, respectively, by the
conditions 
\begin{equation}
\iota _{l}(a)I_{k}=\iota _{k}(a),\quad \pi _{l}(B_{k})=\pi
_{k}(B)E_{l},\quad \forall l\in L,\quad k\subseteq l,
\end{equation}%
where $a\in \mathfrak{A}_{l}$, $B_{k}=B\times \mathrm{E}_{l\backslash k}$, $%
B\in \mathcal{B}_{k}$, $E_{l}=\underset{k\subseteq l}{\wedge }P_{k}$. By
virtue of condition 4 of Definition 2, the representations $\iota _{l}$ and $%
\pi _{l}$ have a common nondecreasing unit $\iota _{l}(1)=E_{l}=\pi _{l}(%
\mathrm{E}_{l})$, $E_{l}\leq E_{l^{\prime }}$, if $l\leq l^{\prime }$ and
are self-consistent, 
\begin{equation}
\iota _{l}(a)=\iota _{l^{\prime }}(a)E_{l},\quad \pi _{l}(B_{k})=\pi
_{l^{\prime }}(B_{k}^{\prime })E_{l},\quad l<l^{\prime }\in L,
\end{equation}%
where $a\in \mathfrak{A}_{l^{\prime }}$, $B_{k}^{\prime }=B\times \mathrm{E}%
_{l^{\prime }\backslash k}$, $B\in \mathcal{B}_{k}$, $k\subseteq l\cup
l^{\prime }$. In accordance with condition 4 of Definition 2 they satisfy
the commutation condition $\pi _{l}(\mathcal{B}_{l})\subseteq \iota _{l}(%
\mathfrak{A}_{l})^{\prime }$ for every $l\in L$, determining uniquely
through the relations (4.1) the original $\mathcal{K}$-process if it
satisfies the conditions 
\begin{equation}
\underset{l\supseteq k}{\vee }E_{l}=P_{k},\qquad \underset{l\supseteq k}{%
\wedge }E_{l}=I_{k},\qquad l\in L,\quad k\in K,
\end{equation}%
and if $\mathfrak{A}_{k}=\underset{l\supseteq k}{\vee }\mathfrak{{A}_{l}}$.
The latter condition in accordance with Theorem 3 is sufficient for the
existence of the reconstruction of the process possessing the properties
(4.3) for $k\neq \emptyset $; the necessary and sufficient condition for
(4.3) to hold for $k=\emptyset $ is the condition (3.2) of asymptotic
noncorrelation with the distant past relative to $\mathcal{K}$. It follows
from the inequality 
\begin{equation}
(\zeta b|\left( E_{l}-\mathrm{P}\right) |\zeta b)\geq \left\vert \sum (\zeta
_{i}|\kappa (b^{i},b)-\kappa (b^{i},e)\kappa (e,b)|\zeta )\right\vert ^{2},
\end{equation}%
as $E_{l}\rightarrow \mathrm{P}$ if $l\rightarrow \emptyset $, where $\zeta
_{i}\in \mathcal{K}$, $b^{i}\in \mathfrak{B}_{l}$ and $\sum (\zeta
_{i}|\kappa (b^{i},b^{i^{\prime }})\zeta _{i^{\prime }})\leq 1$. This
inequality is obtained in the the proof (Sec. 3.1) of the sufficiency of the
condition (3.2) if one uses only the property $E_{l}|\xi _{i}b^{i})=|\xi
_{i}b^{i})$ equivalent to the condition $E_{l}F(b)=F(b)$ for all $b\in 
\mathfrak{B}_{l}$, which is obvious in the representation (2.7) for the
"Feynman integral" $F(\hat{\mathbf{b}})=F_{\emptyset }^{\Lambda }(\mathbf{b}%
) $, $\Lambda \leq l$. Note that if the semigroup $S$ acts on $T$
quasitransitively in the sense that for each $\mathrm{t}\geq \mathrm{t}%
^{\prime }$ there exists $s\in S$ such that $s\mathrm{t}\leq \mathrm{t}%
^{\prime }$, and the conditions of $S$-covariance 
\begin{equation}
V_{s}\iota _{l}(a^{s})=\iota _{sl}(a)V_{s}E_{l},\qquad V_{s}\pi
_{l}(a^{s})=\pi _{sl}(a)V_{s}E_{l},
\end{equation}%
following from the local conditions (2.5) and condition 5 of Definition 2,
are satisfied, the regularity $\underset{l\in L}{\wedge }E_{l}=\mathrm{P}$
ensures generation of the whole space $\mathcal{H}$ by the action of the
operators $U^{s}=V_{s}^{\ast }$, $s\in S$, on any "neighborhood" $\mathcal{H}%
_{l}=E_{l}\mathcal{H}$ of the "infinitely distant initial subspace" $%
\mathcal{K}=\mathrm{P}\mathcal{H}$. In the scalar case $\mathrm{P}=\mathrm{P}%
_{\xi }$ ($\mathcal{K}\simeq \mathbb{C}$) the condition $\underset{l\in L}{%
\wedge }E_{l}=\mathrm{P}_{\xi }$ is a noncommutative analog of the
Kolmogorov regularity condition (law "0 or 1") of a classical stochastic
process in the sense \cite{12}.

\begin{definition}
An $\mathcal{H}$-process $(\mathcal{B},\pi )$ with respect to $(\mathfrak{A}%
,\iota )$, is said to be canonical if $\underset{l\supseteq k}{\vee }%
E_{l}=P_{k}$, $k\in K$, and regular if $\underset{l\supseteq k}{\wedge }%
E_{l}=I_{k}$, $k\in K$, where 
\begin{equation*}
\underset{k\subseteq l}{\wedge }P_{k}=E_{l}=\underset{k\subseteq l}{\vee }%
I_{k},\quad l\in L,\quad P_{\emptyset }=\mathrm{I},\quad I_{\emptyset }=%
\mathrm{P}.
\end{equation*}
\end{definition}

2. A model of a $\mathcal{H}$-process $(\mathcal{B},\pi )$ with respect to $(%
\mathfrak{A},\iota )$ is any $\mathcal{H}^{1}$-process $(\mathcal{B},\pi
^{1})$, with respect to $(\mathfrak{A},\iota ^{1})$ on some Hilbert space $%
\mathcal{H}^{1}\supseteq \mathcal{K}$ for which there exists an isometry $U:%
\mathcal{H}^{1}\rightarrow \mathcal{H}$, $U^{\ast }U=I^{1}$, intertwining
the corresponding representations: 
\begin{equation}
U\iota _{k}^{1}(a)=\iota _{k}(a)UI_{k}^{1},\quad U\pi _{k}^{1}(B)=\pi
_{k}(B)UP_{k}^{1},\quad k\in K.
\end{equation}%
Note that an arbitrary model of a canonical or regular $S$-covariant process
is in general neither canonical, nor regular, nor $S$-covariant; however, it
satisfies the condition $\underset{l\in L}{\wedge }E_{l}=\mathrm{P}$ of
`regularity at the origin" $k=\emptyset $, and will be $S$-covariant if the
condition $UV_{s}^{1}=V_{s}U$ determines for every $s\in S$ an operator $%
V_{s}^{1}:\mathcal{H}^{1}\rightarrow \mathcal{H}^{1}$ for which 
\begin{equation}
V_{s}^{1}I_{k}^{1}=I_{sk}^{1}V_{s}^{1}I_{k}^{1},\quad
V_{s}^{1}P_{k}^{1}=P_{sk}^{1}V_{s}^{1}P_{k}^{1},\quad k\in K.
\end{equation}%
This always holds for a unitarily equivalent model determined by the
condition of its being able to be modeled by the original process with
respect to the adjoint operator $U^{\ast }$ as isometry $\mathcal{H}%
\rightarrow \mathcal{H}^{1}$, i.e., for $I_{k}=UI_{k}^{1}U^{\ast }$ and $%
P_{k}=UP_{k}^{1}U^{\ast }$ for all $k\in K$.

By a subprocess of a $\mathcal{H}$-process $(\mathcal{B},\pi )$ with respect
to $(\mathfrak{A},\iota )$ we understand any model of it determined by
subrepresentations 
\begin{equation}
\iota _{k}^{1}(a)=\iota _{k}(a)I_{k}^{1},\quad \pi _{k}^{1}(B)=\pi
_{k}(B)P_{k}^{1},\quad k\in K,
\end{equation}%
on some subspace $\mathcal{H}^{1}\subseteq \mathcal{H}$ ($U$ is injection,
e.g. $U=\mathrm{I}$, if $\mathcal{H}^{1}=\mathcal{H}$). It necessary for the
subspace $\mathcal{H}^{1}$ to be invariant with respect to $V_{s}$, $s\in S$%
, if the subprocess of the $S$-covariant process is to remain $S$-covariant.
This invariance holds in particular for the minimal subprocess determined by
the least of the orthogonal projectors $I_{k}^{1}$, $P_{k}^{1}$ satisfying
the conditions $I_{k}^{1}F(b)=F(b)$, $b\in \mathfrak{B}_{k}$ and $%
P_{k}^{1}F(b)=F(b)$, $b\in \mathfrak{B}_{k[}$ (this is a canonical process,
called the \emph{minimal modification}), and for the minimal regular
subprocess determined for every regular process by the orthogonal projectors
(4.3) as canonical regular process with respect to the least of the
orthogonal projectors $E_{l}^{1}$ that satisfy the conditions $%
E_{l}^{1}F(b)=F(b)$, $b\in \mathfrak{B}_{l}$, $l\in L$. Note that the
(minimal) subrepresentations $\iota _{l}^{\ast }$ are normal and can be
uniquely extended to von Neumann algebras $\mathfrak{A}_{l}^{\prime \prime }$%
, and similarly the representations $\pi _{l}^{1}$ can be uniquely extended
to $\sigma $-representations of the $\sigma $-algebras $\sigma (\mathcal{B}%
_{l})$, which are generated on $\mathrm{E}_{l}$ by the semirings $\mathcal{B}%
_{l}$ by $\sigma $-additive extension of the original $\sigma $%
-representations $\pi _{l}$, $l\in L$ to $\sigma (\mathcal{B}_{l})$.

The physical equivalence, according to which two quantum stochastic
processes are equivalent if they are statistically indistinguishable in due
course of all possible successive measurements performed on the quantum
model in the causal order is much weaker than the strong unitary
equivalence. This leads to the concept of \emph{wide equivalence} of quantum
stochastic processes described by the same families $\kappa =\{\kappa
^{\Lambda }\}$ of correlation kernels, i.e., coinciding as the processes in
the wide sense. Two processes $\pi $ and $\pi ^{1}$ over the same family $%
\mathcal{B}$ that are represented with respect to the systems $(\mathfrak{A}%
,\iota )$ and $(\mathfrak{A},\iota ^{1})$, respectively, on $\mathcal{H}$
and $\mathcal{H}^{1}$ are said to be equivalent in the wide sense with
respect to the family $\mathfrak{A}$ if they are described by the identical
multikernels as operator-valued functions on the same initial Hilbert space $%
\mathcal{K}$, i.e. if $\kappa (b,b^{\prime })=\kappa ^{1}(b,b^{\prime })$
for all $b,b^{\prime }\in \mathfrak{B}$.

3. The relationship between equivalences in the narrow and wide senses is
established by the following theorem, from which follows the uniqueness up
to unitary equivalence of the minimal (regular) process determining the
canonical reconstruction of Sec. 3.

\begin{theorem}
Let $(\mathcal{B},\pi )$ be a $\mathcal{H}$-process with respect to the
system $(\mathfrak{A},\iota )$ and $(\mathcal{B},\pi ^{1})$ be its isometric
model on $\mathcal{H}^{1}\supseteq \mathcal{H}$ with respect to $(\mathfrak{A%
},\iota ^{1})$. Then these processes are equivalent in the wide sense.
Conversely: The minimal (regular) modification of the $\mathcal{H}$-process $%
(\mathcal{B},\pi )$ with respect to $(\mathfrak{A},\iota )$ is an isometric
model of any other equivalent in the wide sense (regular) $\mathcal{H}$%
-process $(\mathcal{B},\pi ^{1})$ with respect to $(\mathfrak{A},\iota ^{1})$
and is unitarily equivalent to this process if the latter is minimal
(regular).
\end{theorem}

\begin{proof}
Let $U$ be the isometry $\mathcal{H}^{1}\rightarrow \mathcal{H}$ that
determines the modeling relation (4.6), and $F^{1}(b)=\pi
_{k_{n}}^{1}(B_{n})\cdots \pi _{k_{1}}^{1}(B_{1})\mathrm{P}^{1}$ be the
`Feynman integral" that determines the decomposition $\kappa
^{1}(b,b^{\prime })=F^{1}(b)^{\ast }F^{1}(b^{\prime })$. Representing $%
F^{1}(b)$ in the form (2.7), and taking into account (4.6), we obtain $%
UF^{1}(b)=\pi _{k_{n}}(B_{n})\cdots \pi _{k_{1}}(B_{1})U\mathrm{P}^{1}$, and
since for $k\neq \emptyset $ $Ua=aU\mathrm{P}^{1}$ for any $a\in \mathfrak{A}%
_{\emptyset }=\mathcal{B}\left( \mathcal{K}\right) $, the isometric operator 
$U$ on the subspace $\mathcal{K}=\mathrm{P}^{1}\mathcal{H}^{1}$ is a
multiple of the unit: $U\mathrm{P}^{1}=c\mathrm{P}^{1}$, where $|c|=1$. It
follows that $UF^{1}(b)=cF(b)$ and $\kappa ^{1}(b,b^{\prime })=\kappa
(b,b^{\prime })$ for all $b,b^{\prime }\in \mathfrak{B}$, i.e., the relation
of isometric modeling of the processes implies their equivalence in the wide
sense.

Conversely: Suppose the processes $\pi $ and $\pi ^{1}$ have the same
kernels $\kappa $ and $\kappa ^{1}$. By virtue of the transitivity of the
modeling relation, it is sufficient to prove unitary equivalence of the
minimal (regular) modification $\pi ^{\ast }$ of the process $\pi $ and the
minimal (regular) process $\pi ^{1}$. Taking into account the isometry 
\begin{equation*}
(F^{1}(b)\zeta \mid F^{1}(b^{\prime })\zeta ^{\prime })=(F^{\ast }(b)\zeta
\mid F^{\ast }(b^{\prime })\zeta ^{\prime }),\quad \forall b,b^{\prime }\in 
\mathfrak{B},\quad \zeta ,\zeta ^{\prime }\in \mathcal{K}
\end{equation*}%
due to $\kappa ^{1}=\kappa $ of the vector systems $F^{1}(b)\zeta $ and $%
F(b)\zeta =F^{\ast }(b)\zeta $ that generate the minimal subspaces $\mathcal{%
H}^{1}$ and $\mathcal{H}^{\ast }\subseteq \mathcal{H}$, which are invariant
with respect to the corresponding processes $\pi ^{1}$ and $\pi ^{\ast }$,
we can extend the one-to-one correspondence $F^{1}(b)\zeta \mapsto F^{\ast
}(b)\zeta $ of these systems to a unitary operator $U:\mathcal{H}%
^{1}\rightarrow \mathcal{H}^{\ast }$. We then obtain for any $B\in \mathcal{B%
}_{l}$, $l\in L$, 
\begin{equation*}
U\pi _{l}^{1}(B)F^{1}(b)\zeta =F(bB)\zeta =\pi _{l}^{\ast }(B)UF^{1}(b)\zeta
,\quad \forall b\in \mathfrak{B}_{l}
\end{equation*}%
and accordingly for any $a\in \mathfrak{A}_{l}$, $l\in L$, 
\begin{equation*}
U\iota _{l}^{1}(a)F^{1}(b)\zeta =F(b)a\zeta =\iota _{l}^{\ast
}(a)UF^{1}(b)\zeta ,\quad \forall b\in \mathfrak{B}_{l},
\end{equation*}%
from which there follows the unitary equivalence of the systems of
representations $\pi ^{1}$ and $\pi ^{\ast }$ and, respectively, $\iota ^{1}$
and $\iota ^{\ast }$. The proof is completed.
\end{proof}

4. We shall say that a $\mathcal{H}$-process is\/ conditionally Markov, or
simply \emph{dynamic} if it has the property 
\begin{equation}
E_{l}\pi _{l^{\prime }}(\mathcal{B}_{l^{\prime }})E_{l}\subseteq \pi _{l}(%
\mathcal{B}_{l})\vee \iota _{l}(\mathfrak{A}_{l}),\quad \forall l\leq
l^{\prime }\in L,
\end{equation}%
which is obviously inherited by any isometric model of it.

Note that from (4.9) we immediately deduce the regression property (3.3)
with respect to the W*-algebras $\mathcal{A}_{l}=\pi _{l}(\mathcal{B}%
_{l})\vee \iota _{l}(\mathfrak{A}_{l})$ on $\mathcal{H}_{l}=E_{l}\mathcal{H}$
for $\rho _{l}(A)=\mathrm{P}A\mathrm{P}$, $A\in \mathcal{A}_{l}$, and $%
\theta _{ll^{\prime }}(A)=E_{l}AE_{l}$, $A\in \mathcal{A}_{l^{\prime }}$. By
virtue of the fundamental theorem, Theorem 3, the regression condition (3.3)
is also sufficient for existence of the dynamic representation (4.9), and
the minimal canonical realization reconstructing the process satisfying this
weak Markovianity condition is dynamic, i.e. conditionally Markov in the
strong sense. The dynamic mappings $\rho _{l}:\mathcal{A}_{l}\rightarrow 
\mathcal{B}\left( \mathcal{K}\right) $, $\theta _{ll^{\prime }}:\mathcal{A}%
_{l^{\prime }}\rightarrow \mathcal{A}_{l}$ obviously satisfy the conditions
of self-consistency 
\begin{equation*}
\rho _{l}\circ \theta _{ll^{\prime }}=\rho _{l^{\prime }},\;\;\;\;\;\theta
_{l_{0}l}\circ \theta _{ll^{\prime }}=\theta _{l_{0}l^{\prime }}
\end{equation*}%
for $l_{0}\leq l\leq l^{\prime }$ and are $S$-covariant%
\begin{equation*}
\rho _{l}(A^{s})=\rho _{sl}(A)^{s},\;\;\;\theta _{ll^{\prime
}}(A^{s})=\theta _{sl,sl^{\prime }}(A)^{s},
\end{equation*}%
where $A^{s}=V_{s}^{\ast }AV_{s}$. The condition of regularity $%
E_{l}\rightarrow \mathrm{P}$ at the origin $l\rightarrow \emptyset $ of the
dynamic process $\left( \rho _{l},\theta _{ll^{\prime }}\right) $ is
obviously equivalent to the relaxation, or mixing property $\theta
_{l_{0}l}(A)\rightarrow \rho _{l}(A)$, $A\in \mathcal{A}_{l}$, $l_{0},l\in L$
for $l_{0}\rightarrow \emptyset $. This allows the dependence of the limit
relaxed states%
\begin{equation*}
\omega _{l}(A)=\underset{l_{0}\rightarrow \emptyset }{\lim }(\eta \mid
\theta _{l_{0}l}(A)\eta )=\left( \xi \mid \rho _{l}\left( A\right) \xi
\right)
\end{equation*}%
on $\mathcal{A}_{l}$ only on the initial states on the controlling
C*-algebra $\mathfrak{A}_{\emptyset }=\mathcal{B}\left( \mathcal{K}\right) $
determined by vectors $\xi =\mathrm{P}\eta \in \mathcal{K}$ for each $\eta
\in \mathcal{H}$ but not on the initial states on $\mathcal{B}_{l_{0}}$.

An example of a dynamic $\mathcal{H}$-process gives the model of a quantum
coherent filter \cite{17} describing the irreversible process of an optimal
indirect measurements (filtration) of the amplitude of an open (relaxing)
quantum oscillator in a thermostat. Note that the dynamicity condition (4.9)
leads for quantum processes of observation to their weak commutativity $\nu
_{ll^{\prime }}(\mathcal{B}_{l^{\prime }})\subseteq \pi _{l}(\mathcal{B}%
_{l})^{\prime }$ for $l\leq l^{\prime }$, where $\nu _{ll^{\prime
}}(B)=E_{l}\pi _{l^{\prime }}(B)E_{l}$ and we have noted that $\iota _{l}(%
\mathfrak{A}_{l})\subseteq \pi _{l}(\mathcal{B}_{l})^{\prime }\supseteq \pi
_{l}(\mathcal{B}_{l})$, by virtue of which there is complete factorization $%
\kappa (b,b^{\prime })=\mu (b\cdot b^{\prime })$, where $(b\cdot b^{\prime
})(t)=b(t)b(t^{\prime })$ for all $t\in T$. This means that for quantum
dynamic processes defined in the narrow sense ($E_{l}=\mathrm{I}$) ordinary
commutativity must hold, 
\begin{equation}
\pi _{l^{\prime }}(\mathcal{B}_{l^{\prime }})\subseteq \pi _{l}(\mathcal{B}%
_{l})\vee \iota _{l}(\mathfrak{A}_{l})\subseteq \pi _{l}(\mathcal{B}%
_{l})^{\prime },\quad \forall l\leq l^{\prime }\in L,
\end{equation}%
which leads to the strong compatibility condition of the observed events
corresponding to all $l$ and $l^{\prime }\in L$ by virtue of the
directionality of the set $L$. Such representation corresponds to the
dynamic treatment of quantum measurements in the usual, narrow sense as a
strongly Markov process in a "measurement apparatus", described for every
maximal region of simultaneous observation $l\in L$ by a $\ast $-subalgebra $%
\mathfrak{A}_{l}$, which corresponds to the part of the physical medium
which will be causally connected to the observed output process only in the
"future" in accordance with the condition $\pi _{l}(\mathcal{B}%
_{l})\subseteq \iota _{l^{\prime }}(\mathfrak{A}_{l})^{\prime }$ if $l\leq
l^{\prime }$. This latter explains why the $\ast $-algebras $\mathfrak{A}%
_{l} $ are to be assumed nonincreasing, however their representations in $%
\mathcal{H}$ are essential not on the whole Hilbert space but on the
expanding subspaces $\mathcal{H}_{l}\subseteq \mathcal{H}$ containing only
the state vectors to which the large system can pass from the original state 
$\xi \in \mathcal{K}$ in a result of the successive measurements in the
observed subsystem. Such weaker description of the "apparatus" corresponds
to the admission of the irreversibility that is introduced in the process of 
\emph{indirect} successive quantum measurements of noncommuting observables
which can be realized in even a "larger system" without change of the
"future state of the thermostat". This leads to the weakening (4.9) of the
dynamicity condition (4.10), corresponding to weakening of the normalization
condition $\pi _{l}(\mathrm{E}_{l})=\mathrm{I}$ and the assumption of
non-Hamiltonian (conditionally Markovian) evolution. The condition (4.9),
like (4.10), determines independence of the future statistics from the past
if the current state of the observed subsystem and the controlling
thermostat is known, but, in contrast to (4.10), it also permits the
description of irreversible relaxation processes in quantum systems.

5. In conclusion, we show that Kolmogorov reconstructions of classical
stochastic processes and fields arise in the reconstruction process
described in Sec.3 as the minimal processes corresponding to the neglect of
the causality relation on $T=\mathrm{T}$, or the assumption of only
symmetric causality, i.e. an equivalence on $\mathrm{T}$. Note that every
classical process can be represented as a quantum one by $\sigma $%
-homomorphisms $\pi _{k}:\mathcal{B}_{k}\mapsto \mathcal{P}\left( \mathcal{H}%
\right) $ over $\sigma $-algebras $\mathcal{B}_{k}$ with commuting values in
the commutant of the Gel'fand--Naimark--Segal representation $(\mathcal{H}%
,\iota ,\xi )$ of a $\ast $-algebra $\mathfrak{A}$ with a state $\omega $.
In the case of trivial, or symmetric causality the set $L$ contains a single
element $l=T$ -- the factor-set $T=[\mathrm{T}]$, and the set $\mathfrak{B}$
of functions $b:t\mapsto \mathcal{B}_{t}$ can be identified with the
semiring $\mathcal{B}_{l}=\mathcal{B}_{T}$ of cylindrical subsets on $%
\mathrm{E}_{l}=\mathrm{E}_{T}$. Therefore the functional kernel $\kappa
(b,b^{\prime })=\kappa _{T}(b,b^{\prime })$ is completely determined by an
operator-valued (positive, normalized, and $\sigma $-additive) distribution $%
\mu :\mathcal{B}_{T}\rightarrow \mathcal{P}\left( \mathcal{H}\right) $ as $%
\kappa _{T}(b,b^{\prime })=\mu _{T}(b\cdot b^{\prime })$, commuting with a
certain $\ast $-algebra $\mathfrak{A}_{T}\subseteq \mathcal{B}\left( 
\mathcal{K}\right) $, $\mathfrak{A}_{T}=\mu _{T}\left( \mathcal{B}%
_{T}\right) ^{\prime }$ say. The canonical reconstruction of the process
described in the wide sense by such a distribution $\mu _{T}$ reduces to the
construction of a $\sigma $-representation $\pi _{T}:\mathcal{B}%
_{T}\rightarrow \mathcal{P}\left( \mathcal{H}\right) $ on the Hilbert $%
\mathcal{K}$-hull $\mathcal{H}=L_{\mathcal{K}}^{2}(\mathcal{B}_{T})$ of the
Boolean semiring $\mathcal{B}_{T}=\otimes _{t\in T}\mathcal{B}_{t}$ such
that $\mu _{T}\left( B\right) =\mathrm{P}\pi _{T}\left( B\right) \mathrm{P}$%
. It commutes with the corresponding $\ast $-representation $\iota _{T}:%
\mathfrak{A}_{T}\rightarrow \mathcal{B}\left( \mathcal{H}\right) $ of the $%
\ast $-algebra $\mathfrak{A}_{T}$ and determines the resolution of the
identity $\iota _{T}(1)=\mathrm{I}=\pi _{T}(\mathrm{E}_{T})$. The existence
and uniqueness (up to unitary equivalence) of such a representation,
realizing a commutative process with operator-valued probability
distribution $\mu _{T}$ determined by the self-consistent family $\mu
^{\Lambda }(b)=\kappa ^{\Lambda }(b,b)$, $\Lambda \in \mathcal{F}$, are
direct consequences of Theorems 3 and 4. Extending by $\sigma $-additivity
the distribution $\mu _{T}$ to a $\sigma $-measure $\rho _{T}$ on the
generated $\sigma $-algebra $\mathcal{F}_{T}=\sigma (\mathcal{B}_{T})$ of
the space of trajectories $\mathrm{E}_{T}=\times _{t\in T}\mathrm{E}_{t}$,
we obtain a Kolmogorov canonical process $(\mathrm{E}_{T},\mathcal{F}%
_{T},\rho _{T})$ with operator-valued probability measure $\rho _{T}:%
\mathcal{F}_{T}\rightarrow \iota _{T}(\mathfrak{A}_{T})^{\prime }$. The
possibility of this extension for a self-consistent family of
operator-valued distributions $\mu ^{\Lambda }$ and standard Borel spaces $%
\mathrm{E}_{t}$, $t\in T$ with discrete $T=\mathbb{Z}$ was shown by Benioff
in \cite{18}.

Note that the condition of regularity "at the origin" with respect to the
subspace $\mathcal{K}$ as introduced here leads in the case of a
one-dimensional $\mathcal{K}\simeq \mathbb{C}$ and quasitransitive action of 
$S$ on $T$ to the property of complete mixing and ergodicity of the QSP,
which is analogous to $K$ mixing \cite{19} of classical quasiregular systems 
\cite{20}.

I thank V. V. Anshelevich, R. L. Dobrushin, V. P. Maslov, and Ya. G. Sinai
for discussions and interest in the work.

\end{document}